\documentclass[12pt]{article}

\usepackage{graphicx}
\usepackage{epsfig}
\usepackage{ams}

\newcommand{\tr}{^{\prime}}

\def\b#1{\mbox{\boldmath $#1$}}    % - \b: grassetto in formula
\def\bl#1{\mbox{\footnotesize \boldmath {$#1$}}} % - \bl: grassetto per apice
\def\cg#1{\mbox{${\cal #1}$}}
\def\cgl#1{\mbox{\scriptsize {${\cal #1}$}}}

\renewcommand{\th}{\theta}
\newcommand{\al}{\alpha}
\newcommand{\be}{\beta}
\newcommand{\de}{\delta}

\newcommand{\ga}{\gamma}
\newcommand{\si}{\sigma}
\newcommand{\eps}{\varepsilon}
\newcommand{\Th}{\Theta}

\newtheorem{theo}{Theorem}

\begin{document}

\title{Approximate conditional inference for panel logit models\\
allowing for state dependence and\\ unobserved heterogeneity}
\author{Francesco Bartolucci\footnote{Dipartimento di Economia, Finanza e Statistica,
Universit\`a di Perugia, 06123 Perugia, Italy, {\em e-mail:}
bart@stat.unipg.it}\hspace{2mm} and\hspace{1mm} Valentina
Nigro\footnote{Dipartimento di Studi Economico-Finanziari e Metodi
Quantitativi, Universit\`a  di Roma ``Tor Vergata", Via Columbia
2, 00133 Roma, Italy, {\em e-mail:} Valentina.Nigro@uniroma2.it}}
\date{}
\maketitle

\begin{abstract}
\noindent We show that a dynamic logit model for binary panel data
allowing for state dependence and unobserved heterogeneity may be
accurately approximated by a quadratic exponential model, the
parameters of which have the same interpretation that they have in
the true model. We also show how we can eliminate the parameters for
the unobserved heterogeneity from the approximating model by
conditioning on the {\em total scores}, i.e. sum of the response
variables for any individual in the panel. This allows to construct
an approximate conditional likelihood for the dynamic logit model, by
maximizing which we can estimate the parameters for the covariates
and the state dependence. This estimator is very simple to compute
and, by means of a simulation study, we show that it is competitive
in terms of efficiency with the estimator of Honor\'e \& Kyriazidou
(2000). Finally, we outline the extension of the proposed approach to
the case of more elaborated structures for the state dependence and
to that of categorical response variables with more than two
levels.\vspace*{0.5cm}

\end{abstract}

\newpage

\section{Introduction}
An important issue in the econometric literature is the investigation
of the so-called {\em state dependence}, i.e. how the experience of
an event in the past can influence the occurrence of the same event
in the future (see Heckman, 1981a, 1981b). This phenomenon arises in
many economic applications, such as job decision, investment choice
and brand choice. A correct analysis of this phenomenon should take
into account the unobserved heterogeneity between individuals for
what concerns the propensity to experience a certain outcome in all
periods. The latter gives rise to a spurious state dependence that,
as underlined by Heckman, is important to disentangle from the true
state dependence in the analysis of a panel data set, as it can
determine, for instance, different policy implications.

In the case of binary response variables, panel data are usually
analyzed through a dynamic logit or probit model which includes,
among the explanatory variables, the lagged response variable (true
state dependence) and has an individual-specific intercept
(unobserved heterogeneity); see Hsiao (1986) and Arellano \& Honor\'e
(2001), among others. When the latter is considered as a fixed
parameter, the approach suffers from the so-called incidental
parameter problem (Neyman \& Scott, 1948), which leads to
inconsistent estimates of the structural parameters for the
covariates and the true state dependence. For this reason, the
individual specific intercept is frequently considered as a random
parameter (see, for instance, Hyslop, 1999). This requires the
formulation of a certain distribution for this parameter, the
dependence of which on the covariates has to be suitably modelled. In
this case, the problem of the specification of the initial conditions
of the dynamic panel process also arises and the estimation of the
resulting model usually involves multiple integrals which may be
cumbersome to compute.

When a logit model is assumed, an alternative approach for
eliminating the dependence of the joint distribution of the response
variables on the incidental parameters is by conditioning on suitable
statistics. In particular, when the lagged response variable is
omitted from the model, and therefore true state dependence is not
considered, obvious statistics on which conditioning are the sums of
the response variables at individual level. These are sufficient
statistics for the incidental parameters, which, using a terminology
derived from Rasch (1961), will be referred to as total scores. The
resulting maximum likelihood estimator of the other parameters may be
computed by means of a simple Newton-Raphson algorithm and has
optimal asymptotic properties (see  Andersen, 1970, 1972). A
conditional likelihood approach can also be followed when the assumed
logit model includes the lagged response variable. In particular, by
exploiting an intuition of Chamberlain (1985), Honor\'e \& Kyriazidou
(2000) proposed a weighted conditional likelihood that may be used to
consistently estimate the structural parameters. The statistics on
which conditioning are different from the total scores and are such
that a larger number of response configurations does not contribute
to the likelihood. Moreover, the approach requires the specification
of a suitable kernel function for weighting the response
configuration of any subject on the basis of the covariates.

In this paper, we propose a conditional approach for estimating
the parameters of a dynamic logit model for binary panel data
which is based on the approximation of the model through a
particular {\em quadratic exponential model} (Cox, 1972). This
approximation is found by following a method similar to that
adopted by Cox \& Wermuth (1994) in a different context. The
approximating model is in practice a log-linear model for the
conditional distribution of the response variables given the
initial observation and the covariates. The main effects of this
model depend on the covariates and on an individual-specific
parameter for the unobserved heterogeneity, while the two-way
interaction effects are equal to a common parameter when they are
referred to a pair of consecutive response variables and to 0
otherwise. We show that this interaction parameter has the same
interpretation as in the dynamic logit model in terms of {\em
log-odds ratio}, a measure of association between binary variables
which is well known in the statistical literature on categorical
data analysis (Agresti, 2002, Ch. 8).

An interesting feature of the approximating model is that the
parameters for the unobserved heterogeneity may be eliminated by
conditioning on the total scores. This allows to construct an
approximate conditional likelihood for the dynamic logit model, by
maximizing which we obtain an estimator of the structural
parameters. This estimator is simple to compute as the one used in
absence of state dependence and does not require to formulate a
weighting function as the estimator of Honor\'e \& Kyriazidou
(2000) does. The asymptotic properties of this estimator, when the
approximating model holds, are proved on the basis of standard
inferential results (Newey and McFadden, 1994). Under the true
model, instead, they are studied by means of a simulation study
performed along the same lines as Honor\'e \& Kyriazidou (2000).
These simulations show that the proposed estimator is usually more
efficient than their estimator. This is mainly due to the fact
that our approach is based on a likelihood to which a larger
number of response configurations contribute with respect to the
likelihood on which their estimator is based. We also outline the
extension of the proposed approach to the case in which the logit
model includes a second-order lagged response variable and to that
of categorical response variables with more than two levels.

The paper is organized as follows. In the next section we briefly
review the dynamic logit model for binary panel data and describe
the weighted conditional likelihood approach of Honor\'e \&
Kyriazidou (2000); we consider this as a benchmark approach for
the estimation of the model at issue. The proposed approximating
model is described in Section 3, where its conditional
distribution given the total scores is also derived. The resulting
conditional maximum likelihood estimator is described in Section
4, where the asymptotic properties of this estimator under the
approximating model are also illustrated. The results of the
simulation study are shown in Section 5. Finally, in Section 6 we
outline some possible extensions of the proposed approach and in
Section 7 we draw the main conclusions.

All the algorithms described in this paper have been implemented
in {\sc Matlab} functions which are available at the webpage
\verb"www.stat.unipg.it/"$\sim$\verb"bart".
\section{Dynamic logit models for binary panel data}
In the following, we first review the dynamic logit model for
binary panel data and then we discuss conditional maximum
likelihood estimation of its structural parameters.
\subsection{Basic assumptions}
Let $y_{it}$ be a binary random variable equal to 1 if the subject
$i$ ($i=1,\ldots,n$) in the panel makes a certain choice at time $t$
($t=1,\ldots,T$) and to 0 otherwise; also let $\b x_{it}$ be a
corresponding vector of strictly exogenous covariates of size $k$.
The standard econometric model for variables of this type assumes
that
\begin{equation}
y_{it}=1\{\al_i+\b x_{it}\tr\b\be+y_{i,t-1}\ga+\eps_{it}>0\},\quad
i=1,\ldots,n,\quad t=1,\ldots,T,\label{econ1}
\end{equation}
where $1\{\cdot\}$ is the indicator function, $\al_i$ is a fixed or
random individual-specific parameter, the zero-mean random variables
$\eps_{it}$ represent error terms and the initial observations
$y_{i0}$ are assumed to be exogenous. Moreover, $\b\be$ is a vector
of parameters for the covariates and $\ga$ is a parameter measuring
the state dependence effect. The interest is mostly on the last two.
These will be referred to as {\em structural parameters} and, in the
following, will be jointly denoted by $\b\th = (\b\be\tr,\ga)\tr$.
The parameters $\al_i$ are instead considered as {\em incidental
parameters}, the estimation of which is of minor interest.

The typical assumption when the incidental parameters are treated
as fixed parameters is that the errors terms $\eps_{it}$ are
independent and identically distributed conditionally on the
covariates, and with standard logistic distribution. Therefore,
for any subject $i$, the conditional distribution of $y_{it}$
given $\al_i$, $\b X_i=\pmatrix{\b x_{i1} & \cdots & \b x_{iT}}$
and $y_{i0},\ldots,y_{i,t-1}$ may be expressed as
\begin{eqnarray}
p(y_{it}|\al_i,\b X_i,y_{i0},\ldots,y_{i,t-1}) &=& p(y_{it}|\al_i, \b
x_{it},y_{i,t-1}) =\nonumber\\
&=& \frac{\exp[y_{it}(\al_i+\b
x_{it}\tr\b\be+y_{i,t-1}\ga)]}{1+\exp(\al_i+\b x_{it}\tr\b\be+
y_{i,t-1}\ga)},\quad t=1,\ldots,T.\label{rec_logit}
\end{eqnarray}
This is a dynamic logit formulation which implies the following
conditional distribution of the overall vector of response variables
$\b y_i=(y_{i1},\ldots,y_{iT})$ given $\al_i,\b X_i$ and $y_{i0}$:
\begin{equation}
p(\b y_i|\al_i,\b X_i,y_{i0}) = \frac{\exp(y_{i+}\al_i+\sum_t
y_{it}\b x_{it}\tr\b\be+y_{i\times}\ga)}{\prod_t[1+\exp(\al_i+\b
x_{it}\tr\b\be+y_{i,t-1}\ga)]},\label{joint_logit}
\end{equation}
where $y_{i+} = \sum_t y_{it}$ and $y_{i\times}=\sum_t
y_{i,t-1}y_{it}$, with the product $\prod_t$ and the sum $\sum_t$
ranging over $t=1,\ldots,T$.

For what follows, it is important to note some features of the
dependence structure between the response variables in $\b y_i$,
given $\al_i$, $\b X_i$ and $y_{i0}$, implied by the model above.
First of all we have that, for $t=1,\ldots,T-1$, $y_{it}$ is
conditionally independent of any other response variable given
$y_{i,t-1}$ and $y_{i,t+1}$. Moreover, since for $t=1,\ldots,T$ we
have that
\[
\log\frac{p(y_{it}=0|\al_i,\b
x_{it},y_{i,t-1}=0)p(y_{it}=1|\al_i,\b
x_{it},y_{i,t-1}=1)}{p(y_{it}=0|\al_i,\b
x_{it},y_{i,t-1}=1)p(y_{it}=1|\al_i,\b x_{it},y_{i,t-1}=0)}=
\log\frac{\exp(\al_i+\b x_{it}\tr\b\be+\ga)}{\exp(\al_i+\b
x_{it}\tr\b\be)}=\ga,
\]
the parameter $\ga$ for the state dependence is nothing else than the
log-odds ratio between any pair of variables $(y_{i,t-1},y_{it})$,
conditionally on all the other response variables or marginally with
respect to these variables.
\subsection{Conditional inference}
As mentioned in Section 1, an interesting approach for estimating
the fixed effect model illustrated above is based on the
maximization of the conditional likelihood given suitable
statistics. For the case in which the model includes the lagged
response variable, one of the first authors to deal with this
approach was Chamberlain (1985). In particular, he noticed that
when $T=3$ and the covariates are omitted from the model, so that
\[
p(y_{it}|\al_i,y_{i0},\ldots,y_{i,t-1}) =
p(y_{it}|\al_i,y_{i,t-1}) =
\frac{\exp[y_{it}(\al_i+y_{i,t-1}\ga)]}{1+\exp(\al_i+
y_{i,t-1}\ga)},\quad t=1,\ldots,T,
\]
then $p(\b y_i|\al_i,y_{i0},y_{i1}+y_{i2}=1,y_{i3})$ does not depend
on $\al_i$ for any $y_{i0}$ and $y_{i3}$. On the basis of this
conditional distribution it is therefore possible to construct a
likelihood which depends on the response configurations of only
certain subjects (those for which $y_{i1}+y_{i2}=1$) and which allows
to consistently estimate the parameter $\ga$.

The conditional approach above was extended by Honor\'e \&
Kyriazidou (2000) to the case where, as in (\ref{rec_logit}), the
model includes exogenous covariates. In particular, they noticed
that $p(\b y_i|\al_i,\b X_i,y_{i0},y_{i1}+y_{i2}=1,y_{i3})$ is
independent of $\al_i$ provided that $\b x_{i2}=\b x_{i3}$. When
this happens with positive probability, we can therefore estimate
the structural parameters $\b\th$ by maximizing a conditional
likelihood whose logarithm may be expressed as
\[
\sum_i 1\{y_{i1}+y_{i2}=1\}1\{\b x_{i2}-\b x_{i3}=\b 0\}\log[p(\b
y_i|\al_i,\b X_i,y_{i0},y_{i1}+y_{i2}=1,y_{i3})].
\]
For the case in which $p(\b x_{i2}=\b x_{i3})=0$, which typically
occurs in the presence of continuous covariates, Honor\'e \&
Kyriazidou (2000) proposed to estimate $\b\th$ by maximizing a
weighted conditional likelihood defined as above, with the
exception that $1\{\b x_{i2}-\b x_{i3}=\b 0\}$ is substituted by a
Kernel density function $K(\cdot)$. The logarithm of this
likelihood is
\begin{equation}
\sum_i 1\{y_{i1}+y_{i2}=1\}K\left(\frac{\b x_{i2}-\b
x_{i3}}{\si_n}\right)\log[p(\b y_i|\al_i,\b
X_i,y_{i0},y_{i1}+y_{i2}=1,y_{i3})], \label{hk}
\end{equation}
with the bandwidth $\si_n$ a priori fixed. Note that the weight given
to the response configuration of the subject $i$ decreases with the
distance between $\b x_{i2}$ and $\b x_{i3}$, while a large weight is
given to the response configuration of this subject when $\b x_{i2}$
is close to $\b x_{i3}$ and so the property of independence of $p(\b
y_i|\al_i,\b X_i,y_{i0},y_{i1}+y_{i2}=1,y_{i3})$ from $\al_i$
approximately holds.

Honor\'e \& Kyriazidou (2000) also shown how the weighted
conditional approach may be used in the case of $T>3$. In this
case, the approach is based on a {\em pairwise weighted
likelihood} whose logarithm is given by the sum, for any pair of
response variables $(y_{is},y_{it})$, $1<s<t<T$, of an expression
similar to (\ref{hk}) referred to this pair of variables. They
also dealt with dynamic logit models including more than one
lagged response variables and multinomial logit models for
response variables having more than two levels and suggested a
version of the Manski (1987) conditional maximum score estimator
which does not require to formulate any distribution for the error
terms.

Although the weighted conditional estimator of Honor\'e \& Kyriazidou
(2000) is of great interest, its use requires careful choice of the
kernel function and of its bandwidth. This choice obviously affects
the performance of the estimator. Moreover, since only certain
response configurations are considered (e.g. those for which
$y_{i1}+y_{i2}=1$ and $\b x_{i2}$ near to $\b x_{i3}$ in the binary
case with $T=3$), the {\em actual sample size}, i.e. the number of
response configurations which contribute to the likelihood, is
usually much smaller than the nominal sample size $n$. This may
obviously limit the efficiency of the estimator. Moreover, Honor\'{e}
\& Kyriazidou (2000) referred of some problem of applicability of
their approach in presence of time dummies.
\section{Proposed approximation}
In this section, we introduce a quadratic exponential model for
binary panel data that approximates the dynamic logit model
illustrated above and we discuss its main features in comparison to
the true model.
\subsection{Approximating quadratic exponential model}
Along the same lines followed by Cox \& Wermuth (1994) in a
different context, we first take the logarithm of $p(\b
y_i|\al_i,\b X_i,y_{i0})$ as defined in (\ref{joint_logit}), i.e.
\begin{equation}
\log[p(\b y_i|\al_i,\b X_i,y_{i0})] = y_{i+}\al_i+\sum_t y_{it}\b
x_{it}\tr\b\be+y_{i\times}\ga-\sum_t\log[1+\exp(\al_i+\b
x_{it}\tr\b\be+y_{i,t-1}\ga)].\label{logp}
\end{equation}
We then approximate the component which is not linear in the parameter on the
basis of a first-order Taylor series expansion around $\al_i=0$, $\b\be=\b 0$
and $\ga=0$ obtaining
\begin{equation}
\sum_t\log[1+\exp(\al_i+\b
x_{it}\tr\b\be+y_{i,t-1}\ga)]\approx\sum_t[\log(2)+0.5\al_i+0.5\b
x_{it}\tr\b\be]+0.5y_{i*}\ga, \label{approx}
\end{equation}
with $y_{i*}=\sum_t y_{i,t-1}=y_{i0}+y_{i+}-y_{iT}$.

Note that the first term at rhs of the expression above is constant
with respect to $\b y_i$; therefore, by substituting (\ref{approx})
in (\ref{logp}) and renormalizing the exponential of the resulting
expression we obtain the approximation
\begin{equation}
p(\b y_i|\al_i,\b X_i,y_{i0})\approx p^*(\b y_i|\al_i,\b
X_i,y_{i0})=\frac{\exp(y_{i+}\al_i+\sum_ty_{it}\b x_{it}\tr\b\be-0.5y_{i*}\ga+
y_{i\times}\ga)}{\sum_{\bl z} \exp(z_+\al_i+\sum_tz_t\b x_{it}\tr\b\be
-0.5z_*\ga+ z_\times\ga)},\label{quad}
\end{equation}
where the sum at the denominator ranges over all the binary vectors
$\b z=(z_1,\ldots,z_T)$ of dimension $T$ and $z_+$, $z_*$ and
$z_\times$ are defined in an obvious way with $z_0\equiv y_{i0}$. The
approximating model is therefore a quadratic exponential model for
binary variables (Cox, 1972), in which the main effect for $y_{it}$
is equal to $\al_i+\b x_{it}\tr\b\be-0.5\ga$ when $t=1,\ldots,T-1$
and to $\al_i+\b x_{it}\tr\b\be$ when $t=T$ and the two-way
interaction effect for $(y_{is},y_{it})$ is equal to $\ga$ when
$t=s+1$ and to 0 otherwise.

The above expression closely resembles (\ref{joint_logit}), the
main difference being in the denominator which in (\ref{quad})
does not depend on $\b y_i$ and it is simply a {\em normalizing
constant} that may be denoted by $\mu_{it}$. The strong connection
between the two models is clarified by the following Theorem, the
proof of which is given in Appendix.

\begin{theo}
For $i=1,\ldots,n$, the quadratic exponential model (\ref{quad})
implies that the conditional logit of $y_{it}$, given $\al_i,\b
X_i$ and $y_{i0},\ldots,y_{i,t-1}$, is equal to
\begin{equation}
\log\frac{p^*(y_{it}=1|\al_i,\b
X_i,y_{i0},\ldots,y_{i,t-1})}{p^*(y_{it}=0|\al_i,\b
X_i,y_{i0},\ldots,y_{i,t-1})}=\left\{\begin{array}{ll}\al_i+\b
x_{it}\tr\b\be+y_{i,t-1}\ga+\log{\displaystyle\frac{g_{i,t+1}(1)}{g_{i,t+1}(0)}}-0.5\ga
& \mbox{if }\: t<T\cr \al_i+\b x_{it}\tr\b\be+y_{i,t-1}\ga  &
\mbox{if }\: t = T,
\end{array} \right.\label{logit}
\end{equation}
with $g_{it}(z)$ denoting a function depending on the data only
through $\b x_{i,t+1},\ldots,\b x_{i,T}$ and such that
$\log[g_{it}(1)/g_{it}(0)]\approx 0.5\ga$, $t=2,\ldots,T$, where the
approximation is in the sense defined above.

For $i=1,\ldots,n$, model (\ref{quad}) also implies that:
\begin{enumerate}
\item[(i)] $y_{it}$ is conditional independent of
$y_{i0},\dots,y_{i,t-2}$ given $\al_i$, $\b X_i$, $y_{i0}$ and $y_{i,t-1}$
($t=2,\ldots,T$);
\item[(ii)] $y_{it}$ is conditional independent on
$y_{i0},\ldots,y_{i,t-2},y_{i,t+2},\ldots,y_{iT}$, given $\al_i$, $\b X_i$,
$y_{i0}$ and $y_{i,t-1},y_{i,t+1}$ ($t=2,\ldots,T-1$).
\end{enumerate}
\end{theo}\vspace*{0.75cm}

Note that, for $t=T$, logit (\ref{logit}) has exactly the same
parametrization that it has under the dynamic logit model
(\ref{rec_logit}). When $t<T$, this equivalence holds approximately
since $\log[g_{it}(1)/g_{it}(0)]\approx 0.5\ga$. The above Theorem
also implies that
\[
\log\frac{p^*(y_{it}=1|\al_i,\b
X_i,y_{i,t-1}=1)}{p^*(y_{it}=0|\al_i,\b X_i,y_{i,t-1}=1)}-
\log\frac{p^*(y_{it}=1|\al_i,\b
X_i,y_{i,t-1}=0)}{p^*(y_{it}=0|\al_i,\b X_i,y_{i,t-1}=0)}=\ga,\quad
i=1,\ldots,n,\quad t=1,\ldots,T,
\]
and then, under the approximating model, $\ga$ has the same
interpretation that it has under the true model, i.e. log-odds
ratio between any consecutive pair of response variables,
conditionally on all the other response variables or marginally
with respect to these variables. Moreover, the approximating model
reproduces the same conditional independence relations between the
response variables (see (i) and (ii) above) of the dynamic logit
model.

\subsection{Conditional approximating model}
The main advantage of the above approximating model with respect to
the true one is in the availability of minimal sufficient statistics
for the heterogeneity parameters $\al_i$. These statistics are
$y_{i+}$, $i=1,\ldots,n$, which will be referred to as {\em total
scores}. As we show below, in fact, the conditional distribution of
$\b y_i$ given $\b X_i$, $y_{i0}$ and $y_{i+}$ does not depend on
$\al_i$ for any $i$.

First of all note that, under the approximating model,
\[
p^*(y_{i+}|\al_i,\b X_i,y_{i0}) = \sum_{\bl z:z_+=y_{i+}} p^*(\b
z|\al_i,\b X_i,y_{i0})
=\frac{\exp(y_{i+}\al_i)}{\mu_{it}}\sum_{\bl
z:z_+=y_{i+}}\exp(\sum_tz_t\b
x_{it}\tr\b\be-0.5z_*\ga+z_\times\ga),
\]
where the sum is extended to all the binary vectors $\b z$ such
that $z_+=y_{i+}$. Then, after some algebra, the conditional
distribution at issue becomes
\begin{equation}
p^*(\b y_i|\al_i,\b X_i,y_{i0},y_{i+}) =\frac{p^*(\b y_i|\al_i,\b
X_i,y_{i0})}{p^*(y_{i+}|\al_i,\b
X_i,y_{i0})}=\frac{\exp(\sum_ty_{it}\b
x_{it}\tr\b\be-0.5y_{i*}\ga+ y_{i\times}\ga)}{\sum_{\bl
z:z_+=y_{i+}}\exp(\sum_tz_t\b
x_{it}\tr\b\be-0.5z_*\ga+z_\times\ga)}.\label{pcond}
\end{equation}
The expression above does not depend on $\al_i$ and therefore may
also be denoted by $p^*(\b y_i|\b X_i,y_{i+},y_{i0})$. The same
happens for the elements of $\b\be$ corresponding to covariates
which are time-invariant. To make this more clear, consider that
we can multiply the numerator and the denominator of (\ref{pcond})
by $\exp(y_{i+}\b x_{i1}\tr\b\be)$ and, after rearranging terms,
obtain
\begin{equation}
p^*(\b y_i|\b X_i,y_{i0},y_{i+}) =p^*(\b y_i|\b
D_i,y_{i0},y_{i+})= \frac{\exp(\sum_{t>1}y_{it}\b d_{it}\tr\b\be
-0.5y_{i*}\ga+ y_{i\times}\ga)}{\sum_{\bl
z:z_+=y_{i+}}\exp(\sum_{t>1} z_t\b d_{it}\tr\b\be-0.5z_*\ga+
z_\times\ga)},\label{pcond-diff}
\end{equation}
with $\b d_{it}=\b x_{it}-\b x_{i1}$ and $\b D_i=\pmatrix{\b
d_{i2} & \cdots & \b d_{iT}}$. We consequently assume that $\b\be$
does not include the intercept and parameters for the covariates
which are time-invariant because these parameters are not
identified. The same happens for the approach of Honor\'e \&
Kyriazidou (2000).

In Section 4.1 we will show how the structural parameters in $\b\th$
may be estimated by maximizing a conditional likelihood constructed
on the basis of (\ref{pcond-diff}).
\subsection{Improving the approximation}
The quality of approximation (\ref{quad}) depends on the distance of
the parameters from 0 since it is based on the Taylor series
expansion around $\al_i=0$, $\b\be=\b 0$ and $\ga=0$ which is
reported in (\ref{approx}). Obviously, when one or more of these
parameters are far from 0, the quality of the approximation may
considerably be improved by choosing another point of the parameter
space around which performing the Taylor series expansion.

Consider, in particular, the following expansion around $\al_i=0$,
$\b\be=\bar{\b\be}$ and $\ga=0$:
\[
\sum_t\log[1+\exp(\al_i+\b
x_{it}\tr\b\be+y_{i,t-1}\ga)]\approx\sum_t\log\{1+\exp(\b
x_{it}\tr\bar{\b\be})+q_{it}[\al_i+\b
x_{it}\tr(\b\be-\bar{\b\be})]\}+ \sum_t q_{it}y_{i,t-1}\ga,
\]
where $\bar{\b\be}$ is any fixed value of $\b\be$ and
\begin{equation}
q_{it}=\frac{\exp(\b x_{it}\tr\bar{\b\be})}{1+\exp(\b
x_{it}\tr\bar{\b\be})}.\label{qit}
\end{equation}
The latter is equal to the probability that $y_{it}=1$ when the
parameters are fixed as above. This expansion is equal to a
component independent of $\b y_i$ plus $\sum_tq_{it}y_{i,t-1}\ga$
and so, along the same lines as in Section 3.1, it results in
following approximating model
\begin{equation}
p^\dag(\b y_i|\al_i,\b X_i,y_{i0})=\frac{\exp(y_{i+}\al_i+\sum_ty_{it}\b
x_{it}\tr\b\be-\sum_tq_{it}y_{i,t-1}\ga+ y_{i\times}\ga)}{\sum_{\bl z}
\exp(z_+\al_i+\sum_tz_t\b x_{it}\tr\b\be -\sum_tq_{it}z_{t-1}\ga+
z_\times\ga)}.\label{quad2}
\end{equation}
This is a quadratic exponential model which closely resembles the
initial approximating model (\ref{quad}), also in terms of
dependence structure between the response variables and
interpretation of the parameters, and such that the total score
$y_{i+}$ is still a sufficient statistic for $\al_i$. We in fact
have that
\[
p^\dag(\b y_i|\b X_i,y_{i0},y_{i+})=\frac{\exp(\sum_ty_{it}\b
x_{it}\tr\b\be-\sum_tq_{it}y_{i,t-1}\ga+ y_{i\times}\ga)}{\sum_{\bl
z:z_+=y_{i+}}\exp(\sum_tz_t\b x_{it}\tr\b\be-\sum_tq_{it}z_{t-1}\ga+
z_\times\ga)},
\]
which may also be expressed as
\begin{equation}
p^\dag(\b y_i|\b D_i,y_{i0},y_{i+})= \frac{\exp(\sum_{t>1}y_{it}\b
d_{it}\tr\b\be-\sum_tq_{it}y_{i,t-1}\ga+ y_{i\times}\ga)}{\sum_{\bl
z:z_+=y_{i+}}\exp(\sum_{t>1} z_t\b
d_{it}\tr\b\be-\sum_tq_{it}z_{t-1}\ga+
z_\times\ga)}.\label{pcond-diff2}
\end{equation}
On the basis of this distribution, we develop a conditional
likelihood, by maximizing which we obtain an estimator of $\b\th$
which should be more efficient than that based on the conditional
distribution of the initial approximating model, provided that
$\bar{\b\be}$ is suitably chosen.
%In practice we could estimate
%$\b\be$ on the basis of the conditional likelihood derived from
%(\ref{pcond-diff}), then fix $\bar{\b\be}$ at this estimate and
%compute an improved estimate of the parameters on the basis of
%conditional likelihood derived from (\ref{pcond-diff2}). This
%procedure could be iterated until convergence in $\b\be$. This
%approach will be described in more detail in Section 4.3.
This estimator will be illustrated in Section 4.3.

A natural question that rises at this point is why we still rely on
an expansion around a point of the parameter space at which $\al_i=0$
and $\ga=0$, instead of considering a generic point of type
$\al_i=\bar{\al}_i$, $\b\be=\bar{\b\be}$, $\ga=\bar{\ga}$. The first
reason for doing this is that, since within our approach we do not
estimate the parameters $\al_i$, which are ruled out by conditioning
on the total scores, we have no way to choose the $\bar{\al}_i$'s in
practical applications. We could use another estimation method to do
this, but this would complicate considerably the proposed approach.
Moreover, an expansion around $\ga=\bar{\ga}$ results in a model
that, though rather similar to (\ref{quad2}), has sufficient
statistics for the incidental parameters $\al_i$ which differ from
the total scores. On the other hand, a series of simulations, the
results of which are illustrated in Section 5, have shown that the
estimator of $\b\th$ obtained by maximizing the conditional
likelihood based on (\ref{pcond-diff2}) performs considerably better
than that obtained by maximizing the conditional likelihood based on
(\ref{pcond-diff}). In particular, this estimator have a surprisingly
low bias even though samples are generated from a dynamic logit model
of type (\ref{rec_logit}) in which most of the parameters $\al_i$
and/or $\ga$ are far from 0.
\section{Approximate conditional inference}
On the basis of distribution (\ref{pcond-diff}), we can derive an
approximate conditional likelihood for the dynamic logit model
that, for an observed sample $(\b X_i,y_{i0},\b y_i)$,
$i=1,\ldots,n$, has logarithm
\begin{equation}
\ell^*(\b\th) = \sum_i\log[p^*(\b y_i|\b
D_i,y_{i0},y_{i+})]\label{cond-lk}
\end{equation}
and obviously does not depend on the heterogeneity parameters
$\al_i$. Since $\log[p^*(\b y_i|\b D_i,y_{i0},y_{i+})]$ is always
equal to 0 when $y_{i+}=0$ or $y_{i+}=T$, the response
configurations for which this happens do not contribute to
(\ref{cond-lk}). An equivalent expression for $\ell^*(\b\th)$ is
then
\begin{equation}
\ell^*(\b\th) = \sum_i1\{0<y_{i+}<T\}\log[p^*(\b y_i|\b
D_i,y_{i0},y_{i+})].\label{lk_star2}
\end{equation}
The actual sample size is then smaller than the nominal one, but
it is always larger than that we have in the approach of Honor\'e
\& Kyriazidou (2000), which is based on a log-likelihood of type
(\ref{hk}). With $T=3$, for instance, the response configurations
$\b y_i$ omitted from (\ref{lk_star2}) are $(0,0,0)$ and
$(1,1,1)$, whereas also the response configurations $(0,0,1)$ and
$(1,1,0)$ are omitted from (\ref{hk}).

In the following, we show how it is possible to estimate $\b\th$ by
maximizing $\ell^*(\b\th)$ and we study the properties of the
resulting estimator under the approximating model and then, by
simulation, under the true model.
\subsection{Computing the approximate conditional maximum likelihood estimator}
First of all note that distribution (\ref{pcond-diff}) may be
expressed in the canonical exponential family form as
\[
p^*(\b y_i|\b D_i,y_{i0},y_{i+})=\frac{\exp[\b u(\b D_i,y_{i0},\b
y_i)\tr\b\th]}{C(\b\th,\b D_i,y_{i0},y_{i+})},\qquad C(\b\th,\b
D_i,y_{i0},y_{i+}) = \sum_{\bl z:z_+=y_{i+}}\exp[\b u(\b
D_i,y_{i0},\b z)\tr\b\th],
\]
with $\b u(\b D_i,y_{i0},\b y_i)=(\sum_{t>1}y_{it}\b d_{it}\tr,
-0.5y_{i*}+y_{i\times})\tr$. This implies that
\[
\log[p^*(\b y_i|\b D_i,y_{i0},y_{i+})] = \b u(\b D_i,y_{i0},\b
y_i)\tr\b\th-\log[C(\b\th,\b D_i,y_{i0},y_{i+})]
\]
has first derivative vector and second derivative matrix equal,
respectively, to
\[
\b\nabla_{\bl\th}\log[p^*(\b y_i|\b D_i,y_{i0},y_{i+})]=\b v(\b
D_i,y_{i0},\b y_i)\quad\mbox{and}\quad \b\nabla_{\bl\th\bl\th}
\log[p^*(\b y_i|\b D_i,y_{i0},y_{i+})]=-\b S(\b D_i,y_{i0},y_{i+}),
\]
where $\b v(\b D_i,y_{i0},\b y_i)= \b u(\b D_i,y_{i0},\b y_i)-\b m(\b
D_i,y_{i0},y_{i+})$, and with $\b m(\b D_i,y_{i0},y_{i+})$ and $\b
S(\b D_i,y_{i0},y_{i+})$ denoting, respectively, the conditional
expected value and the conditional variance of $\b u(\b D_i,y_{i0},\b
y_i)$ given $\al_i$, $\b D_i$ and $y_{i+}$ under the approximating
model. These are given by
\begin{eqnarray*}
\b m(\b D_i,y_{i0},y_{i+})&=&\sum_{\bl z:z_+=y_{i+}}p^*(\b
z|\b D_i,y_{i0},y_{i+})\b u(\b D_i,y_{i0},\b z)\\
\b S(\b D_i,y_{i0},y_{i+})&=&\sum_{\bl z:z_+=y_{i+}}p^*(\b z|\b
D_i,y_{i0},y_{i+})\b v(\b D_i,y_{i0},\b z)\b v(\b D,y_{i0},\b z)\tr.
\end{eqnarray*}
Consequently, for the conditional log-likelihood $\ell^*(\b\th)$
defined in (\ref{cond-lk}), we have score vector
\begin{equation}
\b s(\b\th)= \sum_i\b v(\b D_i,y_{i0},\b y_i)\label{score}
\end{equation}
and observed information matrix
\begin{equation}
\b J(\b\th)= \sum_i\b S(\b D_i,y_{i0},y_{i+}). \label{inf}
\end{equation}

Note that $\b J(\b\th)$ is always non-negative definite since it
corresponds to the sum of a series of variance-covariance matrices
and therefore $\ell^*(\b\th)$ is always concave. When the sample
size is large enough, this matrix is almost surely positive
definite (see the proof of Theorem 2). In practical application,
we should therefore find that $\ell^*(\b\th)$ is also strictly
concave and has a unique maximum corresponding to the conditional
maximum likelihood estimate
$\hat{\b\th}=(\hat{\b\be}\tr,\hat{\ga})\tr$. This estimate may be
found by a simple Newton-Raphson algorithm. At the $h$th step,
this algorithm updates the estimate of $\b\th$ at the previous
step, $\b\th^{(h-1)}$, as
\[
\b\th^{(h)} = \b\th^{(h-1)}+\b J(\b\th^{(h-1)})^{-1}\b
s(\b\th^{(h-1)}).
\]
Since we also have that the parameter space $\b\Th$ is equal to
$\mathbb{R}^{k+1}$, this algorithm is very simple to implement and
usually converges in a few steps to $\hat{\b\th}$, regardless of the
starting value $\b\th^{(0)}$.
\subsection{Asymptotic properties under the approximating model}
Suppose that the individuals in the samples are independent of each
other with $\al_i$, $\b X_i$, $y_{i0}$ and $\b y_i$ drawn, for
$i=1,\ldots,n$, from the model
\begin{equation}
f_0(\al,\b X, y_0,\b y)=f_0(\al,\b X,y_0)p^*_0(\b y|\al,\b
X,y_0),\label{true_model}
\end{equation}
where $f_0(\al,\b X,y_0)$ denotes the joint distribution of
heterogeneity effect (which is not observed), covariates and initial
observation and $p^*_0(\b y|\al,\b X,y_0)$ denotes the conditional
distribution of the response variables under the approximating
quadratic exponential model (\ref{quad}) when $\b\th=\b\th_0$, with
$\b\th_0$ denoting the true value of its structural parameters.

Under very mild conditions on the distribution of the covariates, we
have that $\hat{\b\th}$ exists, is a $\sqrt{n}$-consistent estimator
of $\b\th_0$ and has asymptotic Normal distribution as
$n\rightarrow\infty$. These results is stated more precisely in the
following Theorem, where $E_0(\cdot)$ denote the expected value under
the true model (\ref{true_model}). As we show in Appendix, the
Theorem may be proved on the basis of standard asymptotic results
(see, for instance, Newey and McFadden, 1994).

\begin{theo} Assume that the distribution $f_0(\al,\b X,y_0)$ is
such that $E_0(\b D\b D\tr)$ exists and is of full rank, with $\b
D=\pmatrix{\b x_2-\b x_1 & \ldots & \b x_T-\b x_1}$. Then, for $T\geq
2$, we have that:
\begin{itemize}
\item (Existence) $\hat{\b\theta}$ exists with probability
approaching 1 as $n\rightarrow \infty$;
\item (Consistency) $\hat{\b\theta}\stackrel{p}{\rightarrow}\b\theta_0$;
\item (Normality)
$\sqrt{n}(\hat{\b\th}-\b\th_0)\stackrel{d}{\rightarrow}N(\b 0,\b
I_0^{-1})$, with $\b I_0=E_0[\b S(\b D,y_0,y_+)]$.
\end{itemize}
\end{theo}\vspace*{0.75cm}

On the basis of the maximum likelihood estimator $\hat{\b\th}$, we
can consistently estimate the matrix $\b I_0$ as
\[
\hat{\b I} = \frac{1}{n}\b J(\hat{\b\th})=\frac{1}{n}\sum_i\hat{\b
S}(\b D_i,y_{i0},y_{i+}),
\]
where $\hat{\b S}(\b D_i,y_{i0},y_{i+})$ is the variance-covariance
matrix of the $i$th score component, computed under the estimated
model. The standard errors of the elements of $\hat{\b\th}$ are then
estimated by the corresponding diagonal elements of $(n\hat{\b
I})^{-1}$ under squared root. This directly derives from Newey \&
McFadden (1994, Sec. 4.2). Note that $n\hat{\b I}$ is equal to $\b
J(\hat{\b\th})$ and so it is obtained as a by-product from the
Newton-Raphson algorithm described in Section 4.1.

Because of the asymptotic normality of $\hat{\b\th}$, it is also
possible to construct an approximate $(1-\al)$-level confidence
interval for any parameter $\be_h$ in $\b\be$ and for $\ga$ as
follows:
\begin{equation}
\hat{\be}_h\mp
z_{\al/2}se(\hat{\be}_h)\quad\mbox{and}\quad\hat{\ga}\mp
z_{\al/2}se(\hat{\ga}),\label{conf1}
\end{equation}
where $se$ denotes the standard error estimated as above and
$z_{\al/2}$ is the $100(1-\al/2)$th percentile of the standard Normal
distribution.

We must again recall that the results above hold under the
approximating quadratic exponential model. Therefore, these results
hold approximately under the dynamic logit model, with the quality of
the approximation depending on the distance between the two models.
To study more precisely these properties under the logit model, we
performed a simulation study along the same lines as Honor\'e \&
Kyriazidou (2000). The results of simulation study are illustrated in
Section 5.
\subsection{Improved approximate conditional estimator}
Once an estimate $\hat{\b\th}$ of $\b\th$ is obtained by maximizing
the log-likelihood $\ell^*(\b\th)$, an improved estimate may be
obtained by maximizing
\[
\ell^\dag(\b\th) = \sum_i\log[p^\dag(\b y_i|\b D_i,y_{i0},y_{i+})],
\]
with $p^\dag(\b y_i|\b D_i,y_{i0},y_{i+})$ denoting the
approximating distribution derived in (\ref{pcond-diff2}) with
$\bar{\b\be}=\hat{\b\be}$. We expect an improvement since
distribution $p^\dag(\b y_i|\al_i,\b X_i,y_{i0})$ should be a
better approximation of the true distribution $p(\b y_i|\al_i,\b
X_i,y_{i0})$ with respect to $p^*(\b y_i|\al_i,\b X_i,y_{i0})$. We
recall that the main difference between the two approximating
distributions is in the correction factor $0.5y_{i*}\ga$ that in
$p^\dag(\b y_i|\al_i,\b X_i,y_{i0})$, and thus also in $p^\dag(\b
y_i|\b D_i,y_{i0},y_{i+})$, is substituted by $\sum_t
q_{it}y_{i,t-1}\ga$, with $q_{it}$ defined in (\ref{qit}).

Maximization of $\ell^\dag(\b\th)$ with respect to $\b\th$ may be
performed on the basis of the same iterative algorithm outlined at
the end of Section 4.1. The only difference is in the computation of
the score vector and the information matrix which are still defined,
respectively, as in (\ref{score}) and (\ref{inf}), but with $\b u(\b
D_i,y_{i0},\b y_i)=(\sum_{t>1}y_{it}\b d_{it}\tr, -\sum_t
q_{it}y_{i,t-1}+y_{i\times})\tr$. Provided that the sample is large
enough, also $\ell^\dag(\b\th)$ is almost surely a strictly concave
function of $\b\th$. This ensures that, in practical applications,
the iterative algorithm converges very easily to the maximum of this
function.

In the algorithm above, the vector $\bar{\b\be}$ used to compute the
probabilities $p^\dag(\b y_i|\b D_i,y_{i0},y_{i+})$ of the
approximating model is held fixed at any iteration. However, it may
be reasonable to update $\bar{\b\be}$ at any step of the algorithm
with the estimate of $\b\be$ obtained at end of the previous step.
This in practice means that the quantities $q_{it}$ are dynamic and
not fixed. As we observed, also this algorithm usually converges very
quickly. We denote the value of $\b\th$ at convergence by
$\tilde{\b\th} = (\tilde{\b\be}\tr,\tilde{\ga})\tr$. To understand if
$\tilde{\b\th}$ represent a real improvement over $\hat{\b\th}$ as an
estimator of $\b\th$, we compared the two estimators by simulation
(see Section 5). Standard errors for the elements of $\tilde{\b\th}$
may be estimated on the basis $(n\tilde{\b I})^{-1}$, where
$n\tilde{\b I}$ is an estimate of the information matrix at
$\tilde{\b\th}$ which is directly produced by the above iterative
algorithm. From these standard errors it is possible to construct
approximate confidence intervals for $\b\th$ as described in the
previous section, i.e.
\begin{equation}
\tilde{\be}_h\mp
z_{\al/2}se(\tilde{\be}_h)\quad\mbox{and}\quad\tilde{\ga}\mp
z_{\al/2}se(\tilde{\ga}).\label{conf2}
\end{equation}
\section{Simulation study of the proposed estimators}
In this section, we illustrate a simulation study carried out to
assess the finite sample properties of the proposed estimators under
the dynamic logit model (\ref{rec_logit}). In order to give more
comparability to our work with the previous literature, we decided to
follow the same simulation design adopted by Honor\'e \& Kyriazidou
(2000), to whom we refer for a more detailed description of this
design. The results concern both the estimator $\hat{\b\th}$, built
on the basis of the initial approximation and described in Section
4.1 ({\em basic conditional estimator}, for short), and the estimator
$\tilde{\b\th}$, built on the basis of the improved approximation and
illustrated in Section 4.3 ({\em improved conditional estimator}, for
short). These results also concern the confidence intervals that may
be constructed, following (\ref{conf1}) and (\ref{conf2}), based
around these estimators.
\subsection{Benchmark design}
Under the {\em benchmark design} of Honor\'e \& Kyriazidou (2000),
samples of different dimension ($n=250, 500, 1000, 2000, 4000$)
are initially generated from a dynamic logit model for $T=3$ time
occasions, with only one covariate and parameters $\be=1$ and
$\ga=0.5$. The covariate is generated by drawing any $x_{it}$
($i=1,\ldots,n$, $t=0,\ldots,T$) from a Normal distribution with
mean $0$ and variance $\pi^2/3$, while any $\al_i$
($i=1,\ldots,n$) is generated as $(x_{i0}+\sum_t x_{it})/(T+1)$.
To study the sensitivity of the results on $T$ and $\ga$, Honor\'e
\& Kyriazidou (2000) then considered a number of time occasions
$T$ equal to 7 and different values of $\ga$ (0.25, 1, 2).

Within our simulation study, we generated 1000 samples from any of
the models described above and, for each sample, we estimated $\be$
and $\ga$. For both parameters we also constructed a 95\% and a 80\%
confidence interval. The results in terms of {\em mean bias}, {\em
root mean squared error} (RMSE), {\em median bias} and {\em median
absolute error} (MAE) of the estimators are displayed in Table 1 and
2. For any $\ga$, these tables also show the ratio\footnote{It is
computed as the expected proportion of response configuration $\b
y_i$ such that $0<y_{i+}<T$.} between the actual sample size and the
nominal sample size $n$. The results, in terms, of actual coverage
level of the confidence intervals are displayed in Table 3.

\begin{table}\centering
{\footnotesize \caption{Performance of the basic and improved
conditional estimators under some benchmark simulation designs with
$T=3$. Percentual numbers are referred to the ratio between the
actual sample size and the nominal one.}
\begin{tabular}{crrrlrrrrrrrrrrr}\hline\hline
\vspace*{-0.2cm}\\
&& && && \multicolumn4c{Estimation of $\be$}&& \multicolumn4c{Estimation of $\ga$}\\
 \cline{7-10}\cline{12-15}
\vspace*{-0.2cm}\\
&& && && \multicolumn1c{Mean} &  & \multicolumn1c{Median} &
&& \multicolumn1c{Mean} & & \multicolumn1c{Median} \\
\multicolumn1c{$\ga$} & &   \multicolumn1c{$n$}  & & Estimator & &
\multicolumn1c{Bias} & \multicolumn1c{RMSE} & \multicolumn1c{Bias} &
\multicolumn1c{MAE} && \multicolumn1c{Bias} & \multicolumn1c{RMSE} &
\multicolumn1c{Bias} & \multicolumn1c{MAE}\\\hline
\vspace*{-0.2cm}\\
 0.25    &&  250 &&  Basic &&  0.039   &   0.144   &   0.025   &   0.110   &&  -0.033  &   0.374   &   -0.036  &   0.299   \\
(60\%)&&      &&      Improved    &&  0.026   &   0.142   &   0.010   &   0.110   &&  -0.017  &   0.360   &   -0.029  &   0.286   \\
&&      500 &&  Basic &&  0.024   &   0.096   &   0.017   &   0.075   &&  -0.038  &   0.274   &   -0.033  &   0.221   \\
&&      &&       Improved    &&  0.010   &   0.093   &   0.003   &   0.073   &&  -0.013  &   0.265   &   -0.012  &   0.213   \\
&&       1000    &&  Basic &&  0.020   &   0.069   &   0.016   &   0.054   &&  -0.034  &   0.191   &   -0.035  &   0.156   \\
&&      &&        Improved    &&  0.005   &   0.066   &   0.002   &   0.053   &&  -0.007  &   0.183   &   -0.011  &   0.146   \\
&&        2000    &&  Basic &&  0.019   &   0.048   &   0.017   &   0.038   &&  -0.040  &   0.134   &   -0.043  &   0.108   \\
&&      &&       Improved    &&  0.004   &   0.045   &   0.002   &   0.035   &&  -0.012  &   0.125   &   -0.011  &   0.099   \\
&&       4000    &&  Basic &&  0.016   &   0.036   &   0.017   &   0.029   &&  -0.040  &   0.101   &   -0.040  &   0.081   \\
&&      &&      Improved    &&  0.001   &   0.033   &   0.001   &   0.026   &&  -0.011  &   0.090   &   -0.011  &   0.072   \\
\vspace*{-0.2cm}\\
    0.5 &&  250 &&  Basic &&  0.055   &   0.155   &   0.035   &   0.116   &&  -0.067  &   0.390   &   -0.079  &   0.313   \\
(57\%)&&      &&          Improved    &&  0.027   &   0.146   &   0.010   &   0.111   &&  -0.026  &   0.361   &   -0.027  &   0.285   \\
      &&  500 &&  Basic &&  0.036   &   0.102   &   0.030   &   0.079   &&  -0.070  &   0.288   &   -0.064  &   0.233   \\
&&      &&          Improved    &&  0.008   &   0.094   &   0.003   &   0.074   &&  -0.021  &   0.272   &   -0.020  &   0.219   \\
      &&  1000    &&  Basic &&  0.033   &   0.075   &   0.029   &   0.059   &&  -0.069  &   0.208   &   -0.072  &   0.167   \\
&&      &&          Improved    &&  0.005   &   0.066   &   0.002   &   0.053   &&  -0.017  &   0.189   &   -0.017  &   0.148   \\
      &&  2000    &&  Basic &&  0.031   &   0.057   &   0.028   &   0.045   &&  -0.074  &   0.152   &   -0.078  &   0.123   \\
&&      &&          Improved    &&  0.003   &   0.047   &   0.001   &   0.037   &&  -0.020  &   0.130   &   -0.013  &   0.103   \\
      &&  4000    &&  Basic &&  0.028   &   0.043   &   0.028   &   0.035   &&  -0.077  &   0.122   &   -0.077  &   0.099   \\
&&      &&          Improved    &&  0.000   &   0.033   &   0.000   &   0.027   &&  -0.023  &   0.095   &   -0.021  &   0.077   \\
\vspace*{-0.2cm}\\
    1   &&  250 &&  Basic &&  0.081   &   0.179   &   0.062   &   0.134   &&  -0.117  &   0.443   &   -0.120  &   0.352   \\
(52\%)&&      &&          Improved    &&  0.029   &   0.154   &   0.012   &   0.116   &&  -0.035  &   0.405   &   -0.039  &   0.319   \\
      &&  500 &&  Basic &&  0.060   &   0.120   &   0.055   &   0.094   &&  -0.127  &   0.333   &   -0.134  &   0.268   \\
&&      &&          Improved    &&  0.011   &   0.101   &   0.005   &   0.079   &&  -0.038  &   0.294   &   -0.043  &   0.234   \\
      &&  1000    &&  Basic &&  0.053   &   0.090   &   0.050   &   0.072   &&  -0.127  &   0.249   &   -0.132  &   0.202   \\
&&      &&          Improved    &&  0.004   &   0.070   &   0.002   &   0.056   &&  -0.034  &   0.203   &   -0.040  &   0.161   \\
      &&  2000    &&  Basic &&  0.050   &   0.071   &   0.046   &   0.057   &&  -0.137  &   0.200   &   -0.143  &   0.165   \\
&&      &&          Improved    &&  0.001   &   0.048   &   -0.003  &   0.038   &&  -0.042  &   0.146   &   -0.043  &   0.116   \\
      &&  4000    &&  Basic &&  0.047   &   0.059   &   0.046   &   0.050   &&  -0.140  &   0.174   &   -0.143  &   0.149   \\
&&      &&          Improved    &&  -0.002  &   0.034   &   -0.002  &   0.027   &&  -0.044  &   0.110   &   -0.045  &   0.089   \\
\vspace*{-0.2cm}\\
    2   &&  250 &&  Basic &&  0.119   &   0.234   &   0.084   &   0.169   &&  -0.144  &   0.592   &   -0.168  &   0.471   \\
(42\%)&&      &&          Improved    &&  0.040   &   0.185   &   0.015   &   0.139   &&  -0.030  &   0.526   &   -0.056  &   0.419   \\
      &&  500 &&  Basic &&  0.086   &   0.154   &   0.070   &   0.116   &&  -0.196  &   0.423   &   -0.216  &   0.345   \\
&&      &&          Improved    &&  0.014   &   0.119   &   0.000   &   0.092   &&  -0.060  &   0.358   &   -0.078  &   0.286   \\
      &&  1000    &&  Basic &&  0.070   &   0.108   &   0.065   &   0.086   &&  -0.200  &   0.326   &   -0.200  &   0.264   \\
&&      &&          Improved    &&  -0.003  &   0.078   &   -0.008  &   0.062   &&  -0.073  &   0.252   &   -0.083  &   0.200   \\
      &&  2000    &&  Basic &&  0.065   &   0.087   &   0.062   &   0.070   &&  -0.211  &   0.279   &   -0.213  &   0.235   \\
&&      &&          Improved    &&  -0.007  &   0.055   &   -0.009  &   0.044   &&  -0.078  &   0.191   &   -0.081  &   0.154   \\
      &&  4000    &&  Basic &&  0.066   &   0.078   &   0.066   &   0.067   &&  -0.211  &   0.247   &   -0.217  &   0.217   \\
&&      &&          Improved    &&  -0.006  &   0.039   &   -0.006  &   0.031   &&  -0.079  &   0.148   &   -0.079  &   0.120   \\
\hline
\end{tabular}}
\end{table}

\begin{table}\centering
{\footnotesize \caption{Performance of the basic and improved
conditional estimators under some benchmark simulation designs with
$T=7$. Percentual numbers are referred to the ratio between the
actual sample size and the nominal one.}
\begin{tabular}{crrrlrrrrrrrrrrr}\hline\hline
\vspace*{-0.2cm}\\
&& && && \multicolumn4c{Estimation of $\be$}&& \multicolumn4c{Estimation of $\ga$}\\
\cline{7-10}\cline{12-15}
\vspace*{-0.2cm}\\
&& && && \multicolumn1c{Mean} &  & \multicolumn1c{Median} &
&& \multicolumn1c{Mean} & & \multicolumn1c{Median} \\
\multicolumn1c{$\ga$} & &   \multicolumn1c{$n$}  & & Estimator & &
\multicolumn1c{Bias} & \multicolumn1c{RMSE} & \multicolumn1c{Bias} &
\multicolumn1c{MAE} && \multicolumn1c{Bias} & \multicolumn1c{RMSE} &
\multicolumn1c{Bias} & \multicolumn1c{MAE}\\\hline
\vspace*{-0.2cm}\\
     0.25   &&  250 &&  Basic &&  0.011   &   0.060   &   0.008   &   0.047   &&  -0.057  &   0.151   &   -0.058  &   0.120   \\
 (92\%)&&      &&          Improved    &&  0.006   &   0.059   &   0.003   &   0.047   &&  -0.006  &   0.152   &   -0.011  &   0.123   \\
&&          500 &&  Basic &&  0.009   &   0.043   &   0.007   &   0.034   &&  -0.056  &   0.115   &   -0.057  &   0.092   \\
&&      &&          Improved    &&  0.003   &   0.042   &   0.002   &   0.033   &&  -0.002  &   0.110   &   -0.002  &   0.088   \\
&&          1000    &&  Basic &&  0.004   &   0.030   &   0.004   &   0.024   &&  -0.056  &   0.090   &   -0.056  &   0.074   \\
&&      &&          Improved    &&  -0.001  &   0.030   &   -0.001  &   0.024   &&  -0.007  &   0.079   &   -0.007  &   0.062   \\
&&          2000    &&  Basic &&  0.006   &   0.022   &   0.006   &   0.018   &&  -0.057  &   0.075   &   -0.055  &   0.063   \\
&&      &&          Improved    &&  0.001   &   0.021   &   0.000   &   0.017   &&  -0.006  &   0.052   &   -0.006  &   0.042   \\
&&          4000    &&  Basic &&  0.006   &   0.016   &   0.007   &   0.013   &&  -0.056  &   0.065   &   -0.055  &   0.057   \\
&&      &&          Improved    &&  0.000   &   0.015   &   0.001   &   0.012   &&  -0.005  &   0.038   &   -0.006  &   0.031   \\
\vspace*{-0.2cm}\\
     0.5 &&  250 &&  Basic &&  0.014   &   0.063   &   0.009   &   0.049   &&  -0.111  &   0.180   &   -0.113  &   0.147   \\
 (91\%)&&      &&          Improved    &&  0.006   &   0.061   &   0.001   &   0.049   &&  -0.008  &   0.153   &   -0.009  &   0.124   \\
&&          500 &&  Basic &&  0.012   &   0.044   &   0.009   &   0.035   &&  -0.112  &   0.151   &   -0.115  &   0.127   \\
&&      &&          Improved    &&  0.003   &   0.042   &   0.001   &   0.034   &&  -0.007  &   0.111   &   -0.009  &   0.089   \\
&&          1000    &&  Basic &&  0.007   &   0.030   &   0.007   &   0.024   &&  -0.112  &   0.133   &   -0.113  &   0.116   \\
&&      &&          Improved    &&  -0.001  &   0.029   &   -0.001  &   0.024   &&  -0.012  &   0.082   &   -0.013  &   0.066   \\
&&          2000    &&  Basic &&  0.009   &   0.024   &   0.009   &   0.019   &&  -0.112  &   0.122   &   -0.110  &   0.112   \\
&&      &&          Improved    &&  0.000   &   0.022   &   0.000   &   0.017   &&  -0.010  &   0.054   &   -0.011  &   0.043   \\
&&          4000    &&  Basic &&  0.009   &   0.018   &   0.009   &   0.014   &&  -0.111  &   0.116   &   -0.110  &   0.111   \\
&&      &&          Improved    &&  0.001   &   0.015   &   0.001   &   0.012   &&  -0.009  &   0.039   &   -0.010  &   0.032   \\
\vspace*{-0.2cm}\\
      1  &&  250 &&  Basic &&  0.012   &   0.065   &   0.007   &   0.051   &&  -0.220  &   0.264   &   -0.227  &   0.229   \\
(87\%)&&      &&          Improved    &&  0.006   &   0.063   &   0.002   &   0.050   &&  -0.020  &   0.157   &   -0.018  &   0.124   \\
&&          500 &&  Basic &&  0.010   &   0.045   &   0.009   &   0.036   &&  -0.218  &   0.243   &   -0.218  &   0.221   \\
&&      &&          Improved    &&  0.004   &   0.044   &   0.004   &   0.035   &&  -0.015  &   0.120   &   -0.014  &   0.095   \\
&&          1000    &&  Basic &&  0.006   &   0.031   &   0.005   &   0.025   &&  -0.219  &   0.232   &   -0.221  &   0.219   \\
&&      &&          Improved    &&  -0.001  &   0.030   &   -0.001  &   0.024   &&  -0.021  &   0.087   &   -0.021  &   0.069   \\
&&          2000    &&  Basic &&  0.007   &   0.023   &   0.005   &   0.018   &&  -0.218  &   0.224   &   -0.217  &   0.218   \\
&&      &&          Improved    &&  0.000   &   0.022   &   -0.001  &   0.017   &&  -0.018  &   0.059   &   -0.018  &   0.047   \\
&&          4000    &&  Basic &&  0.007   &   0.017   &   0.007   &   0.014   &&  -0.219  &   0.222   &   -0.219  &   0.219   \\
&&      &&          Improved    &&  0.000   &   0.016   &   0.000   &   0.013   &&  -0.021  &   0.045   &   -0.021  &   0.036   \\
\vspace*{-0.2cm}\\
     2  &&  250 &&  Basic &&  -0.017  &   0.072   &   -0.022  &   0.058   &&  -0.423  &   0.456   &   -0.431  &   0.425   \\
(76\%)&&      &&          Improved    &&  0.007   &   0.071   &   0.001   &   0.055   &&  -0.065  &   0.191   &   -0.072  &   0.156   \\
&&          500 &&  Basic &&  -0.020  &   0.052   &   -0.023  &   0.042   &&  -0.421  &   0.439   &   -0.423  &   0.421   \\
&&      &&          Improved    &&  0.003   &   0.049   &   0.001   &   0.039   &&  -0.058  &   0.151   &   -0.060  &   0.122   \\
&&          1000    &&  Basic &&  -0.024  &   0.041   &   -0.024  &   0.034   &&  -0.426  &   0.435   &   -0.426  &   0.426   \\
&&      &&          Improved    &&  -0.001  &   0.035   &   -0.002  &   0.028   &&  -0.064  &   0.116   &   -0.066  &   0.095   \\
&&          2000    &&  Basic &&  -0.024  &   0.034   &   -0.025  &   0.028   &&  -0.425  &   0.430   &   -0.424  &   0.425   \\
&&      &&          Improved    &&  -0.001  &   0.024   &   -0.002  &   0.019   &&  -0.064  &   0.092   &   -0.065  &   0.077   \\
&&          4000    &&  Basic &&  -0.024  &   0.030   &   -0.024  &   0.026   &&  -0.428  &   0.430   &   -0.428  &   0.428   \\
&&      &&          Improved    &&  -0.001  &   0.017   &   -0.001  &   0.014   &&  -0.066  &   0.081   &   -0.066  &   0.069   \\
\hline
\end{tabular}}
\end{table}

\begin{table}\centering
{\footnotesize \caption{Coverage levels of the confidence intervals
based on the basic and improved conditional estimators under some
benchmark simulation designs.}
\begin{tabular}{rrrrlrrrrrrrrrrrr}\hline\hline
\vspace*{-0.2cm}\\
&& && && \multicolumn5c{$T=3$}&& \multicolumn5c{$T=7$}\\
 \cline{7-11}\cline{13-17}
\vspace*{-0.2cm}\\
&& && && \multicolumn2c{Interval for $\be$} &&
\multicolumn2c{Interval for $\ga$} &&
\multicolumn2c{Interval for $\be$} && \multicolumn2c{Interval for $\ga$} \\
 \cline{7-8}\cline{10-11}\cline{13-14}\cline{16-17}
 \vspace*{-0.2cm}\\
\multicolumn1c{$\ga$} & &   \multicolumn1c{$n$}  & &
\multicolumn1c{Method} & & \multicolumn1c{95\%} &
\multicolumn1c{80\%} && \multicolumn1c{95\%} & \multicolumn1c{80\%}
&& \multicolumn1c{95\%} & \multicolumn1c{80\%} &&
\multicolumn1c{95\%} & \multicolumn1c{80\%}\\\hline
\vspace*{-0.2cm}\\
    0.25    &&  250 &&  Basic &&  0.944   &   0.802   &&  0.947   &   0.812   &&  0.944   &   0.801   &&  0.931   &   0.754   \\
&&      &&          Improved    &&  0.950   &   0.808   &&  0.950   &   0.802   &&  0.949   &   0.804   &&  0.958   &   0.797   \\
&&          500 &&  Basic &&  0.945   &   0.814   &&  0.953   &   0.798   &&  0.944   &   0.788   &&  0.913   &   0.738   \\
&&      &&          Improved    &&  0.955   &   0.823   &&  0.955   &   0.798   &&  0.944   &   0.792   &&  0.952   &   0.791   \\
&&          1000    &&  Basic &&  0.932   &   0.791   &&  0.952   &   0.794   &&  0.949   &   0.800   &&  0.870   &   0.663   \\
&&      &&          Improved    &&  0.950   &   0.807   &&  0.944   &   0.805   &&  0.956   &   0.796   &&  0.953   &   0.800   \\
&&          2000    &&  Basic &&  0.920   &   0.765   &&  0.945   &   0.765   &&  0.940   &   0.789   &&  0.769   &   0.548   \\
&&      &&          Improved    &&  0.946   &   0.809   &&  0.956   &   0.782   &&  0.948   &   0.804   &&  0.955   &   0.797   \\
&&          4000    &&  Basic &&  0.939   &   0.751   &&  0.936   &   0.754   &&  0.932   &   0.763   &&  0.627   &   0.380   \\
&&      &&          Improved    &&  0.955   &   0.798   &&  0.950   &   0.797   &&  0.955   &   0.790   &&  0.947   &   0.810   \\
\vspace*{-0.2cm}\\
    0.5 &&  250 &&  Basic &&  0.928   &   0.802   &&  0.948   &   0.798   &&  0.940   &   0.798   &&  0.883   &   0.658   \\
&&      &&          Improved    &&  0.946   &   0.825   &&  0.943   &   0.804   &&  0.953   &   0.802   &&  0.951   &   0.813   \\
&&          500 &&  Basic &&  0.937   &   0.811   &&  0.951   &   0.793   &&  0.934   &   0.783   &&  0.813   &   0.564   \\
&&      &&          Improved    &&  0.952   &   0.814   &&  0.955   &   0.793   &&  0.941   &   0.800   &&  0.949   &   0.801   \\
&&          1000    &&  Basic &&  0.916   &   0.758   &&  0.946   &   0.769   &&  0.953   &   0.796   &&  0.658   &   0.396   \\
&&      &&          Improved    &&  0.953   &   0.787   &&  0.951   &   0.810   &&  0.953   &   0.802   &&  0.946   &   0.787   \\
&&          2000    &&  Basic &&  0.896   &   0.734   &&  0.913   &   0.745   &&  0.928   &   0.778   &&  0.380   &   0.159   \\
&&      &&          Improved    &&  0.952   &   0.799   &&  0.950   &   0.781   &&  0.945   &   0.801   &&  0.952   &   0.810   \\
&&          4000    &&  Basic &&  0.878   &   0.666   &&  0.867   &   0.669   &&  0.918   &   0.720   &&  0.131   &   0.029   \\
&&      &&          Improved    &&  0.959   &   0.798   &&  0.941   &   0.782   &&  0.950   &   0.802   &&  0.948   &   0.797   \\
 \vspace*{-0.2cm}\\
    1   &&  250 &&  Basic &&  0.919   &   0.796   &&  0.941   &   0.792   &&  0.941   &   0.800   &&  0.675   &   0.412   \\
&&      &&          Improved    &&  0.946   &   0.827   &&  0.949   &   0.799   &&  0.945   &   0.794   &&  0.947   &   0.798   \\
&&          500 &&  Basic &&  0.912   &   0.763   &&  0.937   &   0.763   &&  0.938   &   0.800   &&  0.479   &   0.226   \\
&&      &&          Improved    &&  0.946   &   0.811   &&  0.948   &   0.793   &&  0.946   &   0.807   &&  0.950   &   0.791   \\
&&          1000    &&  Basic &&  0.875   &   0.711   &&  0.913   &   0.721   &&  0.940   &   0.815   &&  0.181   &   0.052   \\
&&      &&          Improved    &&  0.961   &   0.793   &&  0.949   &   0.802   &&  0.945   &   0.808   &&  0.945   &   0.789   \\
&&          2000    &&  Basic &&  0.833   &   0.629   &&  0.847   &   0.627   &&  0.941   &   0.787   &&  0.009   &   0.002   \\
&&      &&          Improved    &&  0.949   &   0.819   &&  0.939   &   0.791   &&  0.952   &   0.799   &&  0.936   &   0.771   \\
&&          4000    &&  Basic &&  0.746   &   0.485   &&  0.729   &   0.469   &&  0.938   &   0.758   &&  0.000   &   0.000   \\
&&      &&          Improved    &&  0.955   &   0.788   &&  0.928   &   0.757   &&  0.948   &   0.808   &&  0.921   &   0.751   \\
\vspace*{-0.2cm}\\
    2   &&  250 &&  Basic &&  0.903   &   0.785   &&  0.948   &   0.788   &&  0.940   &   0.794   &&  0.286   &   0.118   \\
&&      &&          Improved    &&  0.944   &   0.830   &&  0.947   &   0.825   &&  0.940   &   0.833   &&  0.941   &   0.755   \\
&&          500 &&  Basic &&  0.892   &   0.752   &&  0.921   &   0.737   &&  0.926   &   0.771   &&  0.090   &   0.021   \\
&&      &&          Improved    &&  0.952   &   0.830   &&  0.946   &   0.805   &&  0.946   &   0.808   &&  0.931   &   0.749   \\
&&          1000    &&  Basic &&  0.855   &   0.697   &&  0.891   &   0.675   &&  0.879   &   0.687   &&  0.004   &   0.001   \\
&&      &&          Improved    &&  0.956   &   0.797   &&  0.937   &   0.787   &&  0.946   &   0.799   &&  0.896   &   0.703   \\
&&          2000    &&  Basic &&  0.796   &   0.607   &&  0.790   &   0.542   &&  0.824   &   0.592   &&  0.000   &   0.000   \\
&&      &&          Improved    &&  0.961   &   0.785   &&  0.929   &   0.761   &&  0.955   &   0.799   &&  0.840   &   0.616   \\
&&          4000    &&  Basic &&  0.654   &   0.380   &&  0.625   &   0.357   &&  0.708   &   0.441   &&  0.000   &   0.000   \\
&&      &&          Improved    &&  0.948   &   0.781   &&  0.902   &   0.697   &&  0.948   &   0.798   &&  0.686   &   0.445   \\
\hline
\end{tabular}}
\end{table}

For what concerns the bias of the basic estimator $\hat{\be}$, from
Tables 1 and 2 we can see that this bias is always moderate when
$T=3$ and is negligible when $T=7$. For what concerns the efficiency
of $\hat{\be}$, we can note that both RMSE and MAE of this estimator
decrease as $n$ and $T$ grow. In particular they decrease with $n$ at
a rate close to $\sqrt{n}$ and much faster with $T$. This depends on
the fact that the number of observations that contribute to the
approximate conditional likelihood increases more than proportionally
with $T$ because an increase of $T$ also determines and increase of
the actual sample size. Moreover, both RMSE and MAE increase with
$\ga$. This is mainly due to the fact that an increase of $\ga$, when
this is positive, implies a reduction of the actual sample size,
while the approximation on which our approach is based becomes less
sharp. A completely different scenario may be seen for the basic
estimator $\hat{\ga}$ which is always downward biased. Its bias is
not negligible in most of the cases under consideration and tends to
increase with $\ga$ and, surprisingly, with $n$ and $T$. The
dependence on $n$ is much stronger for $T=3$ than for $T=7$. This
bias has obviously a negative effect on the efficiency of the
estimator. More precisely, both RMSE and MAE decrease as $n$ grows at
a rate much slower than $\sqrt{n}$, especially when $T$ and $\ga$ are
large. With $T=7$ and $\ga=2$, for instance, the MAE of $\hat{\ga}$
is close to be constant with respect to $n$ and may be larger than
that for the case in which $T=3$ and $\ga=2$.

For what concerns the improved estimators $\tilde{\be}$ and
$\tilde{\ga}$, Tables 1 and 2 show that these estimators perform, in
terms of bias and efficiency, much better than the basic estimators
illustrated above. In particular, $\tilde{\be}$ has a bias which is
always negligible and its gain in terms of efficiency with respect to
$\hat{\be}$ increases with $n$ and $\ga$ and does not seem to be
strongly affected by $T$. With $T=3$, for instance, $\hat{\be}$ and
$\tilde{\be}$ have the same MAE when $n=250$ and $\ga=0.25$, but the
MAE of the first estimator is more than the double than that of the
second estimator when $n=4000$ and $\ga=2$. The advantage of the
improved estimator $\tilde{\ga}$ over the basic estimator $\hat{\ga}$
is also more evident. Even though $\tilde{\ga}$ is downward biased,
its bias is almost always moderate and seems to increase very slowly
with $n$ and $\ga$ and to decrease as $T$ grows. Moreover, both RMSE
and MAE of $\tilde{\ga}$ decrease as $n$ grows at a rate close to
$\sqrt{n}$ and much faster in $T$ and increase with $\ga$. The gain
in the terms of efficiency of $\tilde{\ga}$ over $\hat{\ga}$
increases with $n$, $T$ and $\ga$. When $T=3$, $n=250$ and
$\ga=0.25$, for instance, the median bias and the MAE of
$\tilde{\ga}$ are equal respectively to -0.029 and 0.286 whereas, for
$\hat{\ga}$, they are equal respectively to -0.036 and 0.299. When
$T=7$, $n=4000$ and $\ga=2$, instead, the median bias and the MAE of
$\tilde{\ga}$ are equal respectively to -0.066 and 0.069, whereas for
$\hat{\ga}$ they are equal respectively to -0.428 and 0.428.

The superiority of the improved estimators over the basic estimators
is confirmed by the behavior of the confidence intervals constructed
around these estimators. In particular, as may be deduced from Table
3, the actual coverage level of the confidence intervals for $\be$
based on $\hat{\be}$ (see (\ref{conf1})) tends to decrease with $n$
and $\ga$ and to increase with $T$. In practice, the actual coverage
level is significantly smaller than the nominal level only when $T=3$
and $\ga\geq 1$. The confidence intervals based on $\tilde{\be}$ (see
(\ref{conf2})) behave even better, with an actual coverage level
which is always very close to the nominal one. Similar conclusions
may be drawn about the confidence intervals for $\ga$. In this case
however, the actual coverage level of the confidence interval based
on $\hat{\ga}$ may be completely inadequate; this is mainly due to
the bias of this estimator. We have a strong improvement with the
confidence intervals based on $\tilde{\ga}$, even though also the
latter may not be width enough when $\ga$ is large. With $T=7$,
$n=1000$ and $\ga=2$, for instance, the 95\% confidence interval
based on $\hat{\ga}$ has a coverage level of 0.004, whereas that of
the confidence interval based on $\tilde{\ga}$ is equal to 0.896.
\subsection{Other designs}
Following Honor\'{e} \& Kyriazidou (2000), we considered other
simulation designs based on the same dynamic logit model used in the
benchmark design with $T=3$, $\ga=0.5$ and $\be=1$. In particular, we
considered the following designs:
\begin{itemize}
\item {\em $\chi^2(1)$ regressor}: the only difference with respect
to the benchmark design is that any $x_{it}$ ($i=1,\ldots,n$,
$t=0,\ldots,T$) is generated from a $\chi^2(1)$ distribution
transformed to have mean 0 and variance $\pi^2/3$;
\item {\em additional regressors}: samples are generated as in the
benchmark design, but three more covariates are used in the
estimation of the parameters. These covariates, which obviously have
no real effect on the response variables, are generated from the same
Normal distribution used to generate $x_{it}$;
\item {\em trending regressors, $T=3$}: the only difference with
respect to the benchmark design is that the covariate is generated as
$x_{it}=\phi(\psi+0.1t+\zeta_{it})$, with $\phi$ and $\psi$ suitably
chosen and where $\zeta_{i0},\ldots,\zeta_{iT}$ follow a Gaussian
AR(1) process with autoregressive coefficient equal to 0.5,
normalized to have variance $\pi^2/3$;
\item {\em trending regressors, $T=7$}: as in the previous design,
but with $T=7$.
\end{itemize}

The results in terms mean bias, RMSE, median bias and MAE are
displayed in Table 4, while the results in terms of actual coverage
level of the confidence intervals are displayed in Table 5. Given
their superiority over the basic estimators, the results concern only
the improved estimators $\tilde{\be}$ and $\tilde{\ga}$ and the
confidence intervals based on these estimators.

\begin{table}[ht]\centering
{\footnotesize \caption{Performance of the improved conditional
estimator under different simulation designs. Percentual numbers are
referred to the ratio between the actual sample size and the nominal
one.}
\begin{tabular}{lrrrrrrrrrrrrrrr}\hline\hline
\vspace*{-0.2cm}\\
&& && \multicolumn4c{Estimation of $\be$}&& \multicolumn4c{Estimation of $\ga$}\\
 \cline{5-8}\cline{10-13}
 \vspace*{-0.2cm}\\
&& && \multicolumn1c{Mean} &  & \multicolumn1c{Median} &
&& \multicolumn1c{Mean} & & \multicolumn1c{Median} \\
\multicolumn1c{Type of design} & &   \multicolumn1c{$n$}  & &
\multicolumn1c{Bias} & \multicolumn1c{RMSE} & \multicolumn1c{Bias} &
\multicolumn1c{MAE} && \multicolumn1c{Bias} & \multicolumn1c{RMSE} &
\multicolumn1c{Bias} & \multicolumn1c{MAE}\\\hline
\vspace*{-0.2cm}\\
regressors $\chi^2(1)$  &&  250 &&  0.020   &   0.157   &   0.006   &   0.123   &&  -0.020  &   0.326   &   -0.026  &   0.261   \\
   (56\%) &&  500 &&  0.007   &   0.106   &   0.002   &   0.084   &&  -0.016  &   0.230   &   -0.017  &   0.184   \\
    &&  1000    &&  0.002   &   0.073   &   -0.002  &   0.058   &&  -0.031  &   0.163   &   -0.028  &   0.130   \\
    &&  2000    &&  -0.001  &   0.052   &   -0.002  &   0.042   &&  -0.024  &   0.113   &   -0.023  &   0.091   \\
    &&  4000    &&  0.000   &   0.039   &   -0.001  &   0.031   &&  -0.024  &   0.080   &   -0.022  &   0.063   \\
\vspace*{-0.2cm}\\
additional regressors   &&  250 &&  0.052   &   0.155   &   0.041   &   0.118   &&  -0.022  &   0.398   &   -0.039  &   0.320   \\
   (57\%)
    &&  500 &&  0.017   &   0.097   &   0.013   &   0.076   &&  -0.015  &   0.257   &   -0.022  &   0.205   \\
    &&  1000    &&  0.013   &   0.064   &   0.013   &   0.051   &&  -0.033  &   0.182   &   -0.037  &   0.147   \\
    &&  2000    &&  0.003   &   0.048   &   0.001   &   0.038   &&  -0.022  &   0.130   &   -0.022  &   0.104   \\
    &&  4000    &&  0.003   &   0.032   &   0.001   &   0.026   &&  -0.016  &   0.090   &   -0.011  &   0.072   \\
\vspace*{-0.2cm}\\
trending regressors,  &&  250 &&  0.030   &   0.171   &   0.016   &   0.129   &&  -0.029  &   0.417   &   -0.036  &   0.328   \\
$T=3$ &&  500 &&  0.013   &   0.117   &   0.001   &   0.092   &&  -0.030  &   0.281   &   -0.028  &   0.225   \\
(42\%)    &&  1000    &&  0.002   &   0.080   &   -0.004  &   0.064   &&  -0.019  &   0.198   &   -0.014  &   0.158   \\
    &&  2000    &&  0.002   &   0.059   &   0.001   &   0.047   &&  -0.034  &   0.145   &   -0.036  &   0.115   \\
    &&  4000    &&  -0.001  &   0.039   &   -0.003  &   0.031   &&  -0.024  &   0.100   &   -0.028  &   0.080   \\
\vspace*{-0.2cm}\\
trending regressors,  &&  250 &&  0.009   &   0.072   &   0.004   &   0.056   &&  -0.015  &   0.168   &   -0.018  &   0.135   \\
$T=7$  &&  500 &&  0.006   &   0.050   &   0.004   &   0.041   &&  -0.013  &   0.122   &   -0.011  &   0.095   \\
(78\%)    &&  1000    &&  0.002   &   0.035   &   0.001   &   0.028   &&  -0.015  &   0.087   &   -0.013  &   0.068   \\
    &&  2000    &&  0.002   &   0.026   &   0.002   &   0.021   &&  -0.014  &   0.060   &   -0.017  &   0.048   \\
    &&  4000    &&  0.002   &   0.018   &   0.001   &   0.015   &&  -0.015  &   0.044   &   -0.015  &   0.036   \\
\hline
\end{tabular}}
\end{table}

\begin{table}[ht]\centering
{\footnotesize \caption{Coverage levels of the confidence intervals
based on the improved conditional estimator under different
simulation designs.}
\begin{tabular}{lrrrcrrrrrrrrrrrr}\hline\hline
\vspace*{-0.2cm}\\
&& && \multicolumn2c{Interval for $\be$} &&
\multicolumn2c{Interval for $\ga$}\\
\cline{5-6}\cline{8-9}
\vspace*{-0.2cm}\\
\multicolumn1c{Type of design} & &   \multicolumn1c{$n$}  & &
\multicolumn1c{95\%} & \multicolumn1c{80\%} && \multicolumn1c{95\%} &
\multicolumn1c{80\%}\\\hline
\vspace*{-0.2cm}\\
regressors $\chi^2(1)$  &&  250 &&  0.947   &   0.815   &&  0.951   &   0.803   \\
    &&  500 &&  0.948   &   0.821   &&  0.948   &   0.798   \\
    &&  1000    &&  0.960   &   0.794   &&  0.940    &   0.802   \\
    &&  2000    &&  0.960   &   0.805   &&  0.947   &   0.803   \\
    &&  4000    &&  0.952   &   0.805   &&  0.934   &   0.779   \\
\vspace*{-0.2cm}\\
additional regressors   &&  250 &&  0.941   &   0.811   &&  0.955   &   0.817   \\
    &&  500 &&  0.942   &   0.800   &&  0.946   &   0.810    \\
    &&  1000    &&  0.945   &   0.803   &&  0.945   &   0.795   \\
    &&  2000    &&  0.950   &   0.816   &&  0.951   &   0.782   \\
    &&  4000    &&  0.946   &   0.794   &&  0.956   &   0.800 \\
\vspace*{-0.2cm}\\
trending regressors,  &&  250 &&  0.951   &   0.826   &&  0.952   &   0.813   \\
$T=3$    &&  500 &&  0.945   &   0.820   &&  0.948   &   0.801   \\
    &&  1000    &&  0.949   &   0.796   &&  0.948   &   0.798   \\
    &&  2000    &&  0.955   &   0.805   &&  0.943   &   0.789   \\
    &&  4000    &&  0.952   &   0.793   &&  0.940    &   0.786   \\
\vspace*{-0.2cm}\\
trending regressors,   &&  250 &&  0.940   &   0.805   &&  0.949   &   0.796   \\
$T=7$    &&  500 &&  0.954   &   0.799   &&  0.945   &   0.815   \\
    &&  1000    &&  0.946   &   0.808   &&  0.942   &   0.800 \\
    &&  2000    &&  0.947   &   0.798   &&  0.945   &   0.801   \\
    &&  4000    &&  0.952   &   0.801   &&  0.941   &   0.785   \\
\hline
\end{tabular}}
\end{table}

On the basis of the results in Table 4 we can conclude that the
improved estimators have not a considerably different behavior with
respect to the benchmark design. Even when the estimators perform
worse, in terms of bias and/or efficiency, with respect to the
benchmark design, the difference is slight. This happens, for the
$\chi^2(1)$ design (limited to $\tilde{\be}$), for the additional
regressors design when $n$ is small and for the trending regressor
design when $T=3$. Occasionally, it also happens that the estimators
perform better with respect to the benchmark design. Limited to
$\tilde{\ga}$, this happens, for instance, for the $\chi^2(1)$
design.

Finally, for what concerns the confidence intervals, we observed that
actual coverage value is always very close to the nominal level for
both parameters $\al$ and $\be$. This confirms the good quality of
the method proposed in Section 4.3 for constructing confidence
intervals, already noticed for the benchmark design.

\subsection{Comparison with the weighted conditional estimator}
An important issue is how the improved version of our approximate
conditional estimator, which we established to be much better than
its basic version, performs in comparison to the weighted conditional
estimator of Honor\'{e} \& Kyriazidou (2000). We then compared their
simulation results with the simulation results illustrated above. An
advantage of our estimator over their estimator, in terms of bias and
efficiency, seems clearly to emerge. The results of this comparison
are summarized in Table 6, which, for certain reference situations
and for both $\be$ and $\ga$, shows the median bias and the MAE of
our estimator in comparison to those of the weighted conditional
estimator. For both estimators, the table also shows the
rate\footnote{For the weighted conditional estimator, this rate is
computed as the expected proportion of pairs of response variables
$(y_{is},y_{it})$, $0<s<t<T$, such that $y_{is}+y_{it}=1$.} between
the actual sample size and the nominal sample size.

\begin{table}\centering
{\footnotesize \caption{Comparison between the weighted and the
improved conditional estimator. Percentual numbers in the first two
columns are referred to actual sample size under the two approaches.
Percentual numbers in the other columns are referred to the reduction
of median bias (in absolute value) and MAE from the first to the
second estimator.}
\begin{tabular}{rrrrrrlrrrrrrrrrr}\hline\hline
\vspace*{-0.2cm}\\
 && && && && \multicolumn2c{Estimation of $\be$} &&
\multicolumn2c{Estimation of $\ga$} \\
 \cline{9-10}\cline{12-13}
 \vspace*{-0.2cm}\\& &     & &   & &  & & \multicolumn1c{ Median} & & &\multicolumn1c{
Median} & \\
\multicolumn1c{$\ga$} & &     \multicolumn1c{$T$}& &
\multicolumn1c{$n$}  & & \multicolumn1c{Estimator} & &
\multicolumn1c{Bias} & \multicolumn1c{MAE} & &\multicolumn1c{Bias} &
\multicolumn1c{MAE}\\
\hline
\vspace*{-0.2cm}\\
  0.5 &&  3   &&  250 &&  Weighted  &&  0.076   &   0.154   &&  -0.039  &   0.403   \\
\multicolumn3c{(37\% - 57\%)}          &&      &&  Approximated    &&  0.010    &   0.111   &&  -0.027  &   0.285   \\
    &&      &&      &&      && (87\%) &   (28\%) &&  (31\%) &   (29\%) \\
    &&      &&  1000    &&  Weighted    &&  0.038   &   0.086   &&  -0.035  &   0.178   \\
    &&      &&      &&  Approximated    &&  0.002   &   0.053   &&  -0.017  &   0.148   \\
    &&      &&      &&      &&  (95\%) &   (38\%) &&  (51\%) &   (17\%) \\
    &&      &&  4000    &&  Weighted    &&  0.019   &   0.044   &&  -0.035  &   0.102   \\
    &&      &&      &&  Approximated    &&  0.000   &   0.027   &&  -0.021  &   0.077   \\
    &&      &&      &&      &&  (100\%)    &   (39\%) &&  (40\%) &   (25\%) \\
\vspace*{-0.2cm}\\
 &&  7   &&  250 &&  Weighted     &&  0.014   &   0.050    &&  -0.053  &   0.131   \\
\multicolumn3c{(43\% - 91\%)}      &&      &&  Approximated     &&  0.001   &   0.049   &&  -0.009  &   0.124   \\
&&      &&      &&      &&  (93\%) &   (2\%)  &&  (83\%) &   (5\%)  \\
&&      &&  1000    &&  Weighted    &&  0.009   &   0.027   &&  -0.041  &   0.075   \\
&&      &&      &&  Approximated    &&  -0.001  &   0.024   &&  -0.013  &   0.066   \\
&&      &&      &&      &&  (89\%)    &   (11\%) &&  (68\%) &   (12\%) \\
&&      &&  4000    &&  Weighted    &&  0.005   &   0.015   &&  -0.033  &   0.039   \\
&&      &&      &&  Approximated    &&  0.001   &   0.012   &&  -0.010   &   0.032   \\
&&      &&      &&      &&  (80\%) &   (20\%) &&  (70\%) &   (18\%) \\
 \vspace*{-0.2cm}\\
2   &&  3   &&  250 &&  Weighted     &&  0.196   &   0.251   &&  -0.056  &   0.620    \\
\multicolumn3c{(26\% - 42\%)}      &&      &&  Approximated     &&  0.015   &   0.139   &&  -0.056  &   0.419   \\
    &&      &&      &&      &&  (92\%) &   (45\%) &&  (0\%)  &   (32\%) \\
    &&      &&  1000    &&  Weighted    &&  0.113   &   0.136   &&  -0.148  &   0.321   \\
    &&      &&      &&  Approximated    &&  -0.008  &   0.062   &&  -0.083  &   0.200 \\
    &&      &&      &&      &&  (93\%)    &   (54\%) &&  (44\%) &   (38\%) \\
    &&      &&  4000    &&  Weighted    &&  0.063   &   0.074   &&  -0.118  &   0.163   \\
    &&      &&      &&  Approximated    &&  -0.006  &   0.031   &&  -0.079  &   0.120    \\
    &&      &&      &&      &&  (90\%)    &   (58\%) &&  (33\%) &   (26\%) \\
\vspace*{-0.2cm}\\
&&  7   &&  250 &&  Weighted     &&  0.016   &   0.064   &&  -0.195  &   0.227   \\
\multicolumn3c{(34\% - 76\%)} &&      &&  Approximated     &&  0.001   &   0.055   &&  -0.072  &   0.156   \\
&&      &&      &&      &&  (94\%) &   (14\%) &&  (63\%) &   (31\%) \\
&&      &&  1000    &&  Weighted    &&  0.016   &   0.034   &&  -0.160   &   0.164   \\
&&      &&      &&  Approximated    &&  -0.002  &   0.028   &&  -0.066  &   0.095   \\
&&      &&      &&      &&  (88\%)    &   (18\%) &&  (59\%) &   (42\%) \\
&&      &&  4000    &&  Weighted    &&  0.006   &   0.017   &&  -0.116  &   0.116   \\
&&      &&      &&  Approximated    &&  -0.001  &   0.014   &&  -0.066  &   0.069   \\
&&      &&      &&      &&  (83\%)    &   (18\%) &&  (43\%) &   (41\%) \\
\hline
\end{tabular}}
\end{table}

From Table 6 we can see that, as regards the parameter $\be$, the
advantage of our estimator $\tilde{\be}$ is particularly evident for
the case $n=250$, $T=3$ and $\ga=2$, case in which $\tilde{\be}$ has
a median bias of 0.015, whereas the weighted conditional estimator
has a median bias of 0.196. For what concerns the efficiency, the
gain of our estimator seems to increase with $n$ and $\ga$ and is
more evident for $T=3$ then for $T=7$. For the case of $n=250$, $T=3$
and $\ga=0.5$, for instance, the reduction of MAE is just of 2\%,
which increases to 58\% for the case in which $n=4000$, $T=3$ and
$\ga=2$. In most of the cases considered in Table 6, the reduction of
MAE is at least of 15\%.

As regards the parameter $\ga$, the reduction of bias is particularly
relevant when $T$ and $\ga$ are large. For instance, with $n=250$,
$T=7$ and $\ga=2$, their estimator has a median bias of -0.195,
whereas our estimator has a median bias of -0.072. Similarly, the
efficiency of our estimator with respect to their estimator seems to
increase with $\ga$, whereas it has not a clear trend in $n$ and $T$.
For instance, with $n=250$, $T=7$ and $\ga=0.5$, the reduction of MAE
from their estimator to our estimator is of 5\%, while it is equal to
41\% for the case of $n=4000$, $T=7$ and $\ga=2$. In most of the
cases considered in Table 6, the reduction of MAE is at least of 25\%
and is usually more evident than for the estimation of $\be$.

The main explanation that we can give for the results above is that,
as may also be deduced from Table 6, the actual sample size used in
our approach is always much larger than that used in the approach of
Honor\'{e} \& Kyriazidou (2000). This difference increases with $\ga$
and $T$. For instance, with $\ga=0.5$ and $T=3$, the actual sample
size used in our approach is about 1.5 times that used in their
approach. This ratio becomes equal to about 2.1 for $\ga=0.5$ and
$T=7$ and to 2.2 for $\ga=2$ and $T=7$. Note however that the gain in
median bias and MAE does not closely follows the gain in the actual
sample size. Other factors have therefore to be taken into
consideration which may affect the performance of the two estimators
in a way that depends on $\ga$ and $T$. We recall, in particular,
that the performance of our estimator depends on the quality of the
approximation we are relying on, while the performance of the
estimator of Honor\'{e} \& Kyriazidou (2000) depends also on the fact
that the response configurations are differently weighted on the
basis of the corresponding covariate configurations and that, for
$T>3$, they are indeed relying on a pairwise likelihood.
\section{Possible extensions}
In the following, we illustrate two possible extensions of the
proposed approach to the case of dynamic logit models including
more than one lagged response variables and to that of multinomial
logit models for categorical response variables with more than two
levels. In both cases, the approximate conditional inference
outlined in the previous sections may be implemented with minor
adjustments.
\subsection{More than one lagged response variables among the regressors}
Sometimes, it may be interesting to know how long is the dynamics of
a certain phenomenon. In our context, to have the possibility to test
for its length it is necessary to use a dynamic logit model with more
than one lagged response variables.

As an illustration consider the case of two lagged response
variables. The model described in Section 2.1 becomes
\begin{eqnarray}
&&p(y_{it}|\al_i,\b X_i,y_{i,-1},\ldots,y_{i,t-1}) =
p(y_{it}|\al_i,\b
x_{it},y_{i,t-2},y_{i,t-1}) = \nonumber\\
&&\hspace*{1cm}= \frac{\exp[y_{it}(\al_i+\b
x_{it}\tr\b\be+y_{i,t-1}\ga_1+y_{i,t-2}\ga_2)]}{1+\exp(\al_i+\b
x_{it}\tr\b\be+ y_{i,t-1}\ga_1+y_{i,t-2}\ga_2)},\quad
i=1,\ldots,n,\quad t=1,\ldots,T,\label{lag2}
\end{eqnarray}
with $\ga_1$ and $\ga_2$ having an obvious interpretation and
$y_{i,-1}$ and $y_{i0}$ assumed to be exogenous. Along the same lines
as in Section 2.1, it is straightforward to write the distribution of
$\b y_i$, given $\al_i$, $\b X_i$, $y_{i,-1}$ and $y_{i0}$, as
\begin{equation}
p(\b y_i|\al_i,\b X_i,y_{i,-1},y_{i0}) =
\frac{\exp(y_{i+}\al_i+\sum_t y_{it}\b x_{it}\tr\b\be+y_{i\times
1}\ga_1+y_{i\times 2}\ga_2)}{\prod_t[1+\exp(\al_i+\b
x_{it}\tr\b\be+y_{i,t-1}\ga_1+y_{i,t-2}\ga_2)]},\label{joint_logit2}
\end{equation}
where $y_{i\times 1}=\sum_t y_{i,t-1}y_{it}$ and $y_{i\times
2}=\sum_t y_{i,t-2}y_{it}$.

In this case, we can approximate the logarithm of the denominator
with a first-order Taylor series expansion around $\al_i=0$,
$\b\be=\b 0$ and $\ga_1=\ga_2=0$ obtaining
\begin{eqnarray*}
&&\sum_t\log[1+\exp(\al_i+\b
x_{it}\tr\b\be+y_{i,t-1}\ga_1+y_{i,t-2}\ga_2)]\approx\\
&&\hspace*{1cm}\approx\sum_t[\log(2)+0.5\al_i+0.5\b
x_{it}\tr\b\be]+0.5\sum_t(y_{i,t-1}\ga_1+y_{i,t-2}\ga_2).
\end{eqnarray*}
Therefore, by substituting the latter into (\ref{joint_logit2})
and after some algebra, we find that $p(\b y_i|\al_i,\b
X_i,y_{i,-1},y_{i0})$ may be approximated with
\begin{eqnarray*}
p^*(\b y_i|\al_i,\b
X_i,y_{i,-1},y_{i0})=\frac{\exp(y_{i+}\al_i+\sum_t y_{it}\b
x_{it}\tr\b\be-0.5y_{i*1}\ga_1-0.5y_{i*2}\ga_2+
y_{i\times1}\ga_1+y_{i\times2}\ga_2)}{\sum_{\bl z}
\exp(z_+\al_i+\sum_t z_t\b x_{it}\tr\b\be
-0.5z_{*1}\ga_1-0.5z_{*2}\ga_2+
z_{\times1}\ga_1+z_{\times2}\ga_2)},
\end{eqnarray*}
where $y_{i*h}=\sum_t y_{i,t-h}$ and $y_{i\times h}=\sum_t
y_{t-h}y_t$, for $h=1,2$, and $z_{*h}$ and $z_{\times h}$ defined in
a similar way, with $z_{-1}\equiv y_{i,-1}$ and $z_0\equiv y_{i0}$.
The approximating model is therefore a quadratic exponential model in
which the main effect parameter for $y_{it}$ is equal to $\al_i+\b
x_{it}\tr\b\be-0.5\ga_1-0.5\ga_2$ when $t=1,\ldots,T-2$, to $\al_i+\b
x_{it}\tr\b\be-0.5\ga_1$ when $t=T-1$ and to $\al_i+\b
x_{it}\tr\b\be$ when $t=T$; moreover, the two-way interaction effect
for $(y_{is},y_{it})$ is equal to $\ga_1$ when $t=s+1$, to $\ga_2$
when $t=s+2$ and to 0 otherwise. The advantage of this model is that
of having a minimal sufficient statistic for $\al_i$ which is again
$y_{i+}$, so that the conditional distribution of $\b y_i$ given $\b
X_i$, $y_{i,-1}$, $y_{i0}$ and $y_{i+}$ does not depend on $\al_i$.
The estimation of the structural parameters follows by maximizing a
likelihood based on this conditional distribution in a way similar to
that outlined in Section 4.1. In a similar way we can also compute
standard errors for these estimates.

In the case outlined above, it may interesting to test the hypothesis
$\ga_2=0$ under which model (\ref{lag2}) specializes into model
(\ref{rec_logit}). In the present approach, this hypothesis may be
tested in the usual way by using the statistic
$\hat{\ga}_2/se(\hat{\ga}_2)$, where $se(\hat{\ga}_2)$ is the
standard error for $\hat{\ga}_2$ estimated as described in Section
4.1. Under the null hypothesis, this statistic should approximately
have a standard Normal distribution.
\subsection{Categorical response variables}
Suppose that any response variable has $M$, instead of 2, possible
levels, from 0 to $M-1$. The standard econometric model assumed in
this case is the dynamic multinomial logit model
\begin{eqnarray*}
p(y_{it}|\al_i,\b X_i,y_{i0},\ldots,y_{i,t-1}) &=&
p(y_{it}|\al_i,\b
x_{it},y_{i,t-1}) = \\
&=&\frac{\exp(\al_{iy_{it}}+\b
x_{it}\tr\b\be_{y_{it}}+\ga_{y_{i,t-1}y_{it}})}{\sum_m\exp(\al_{im}+\b
x_{it}\tr\b\be_m+\ga_{y_{i,t-1}m})},\quad i=1,\ldots,n,\quad
t=1,\ldots,T,
\end{eqnarray*}
where $\al_{i0}=0$ for any $i$, $\b\be_0=\b 0$ and $\ga_{hm}=0$
whenever $h=0$ or $m=0$. It is now convenient to use a dummy
representation for the response variables $y_{it}$ and so let $\b
a_{it}$ be an $(M-1)$-dimensional vector with all elements equal to
0, apart from the element $a_{itm}$, $m=y_{it}-1$, equal to 1 when
$y_{it}>0$. Thus
\[
p(y_{it}|\al_i,\b x_{it},y_{i,t-1}) =
\frac{\exp(\sum_ma_{itm}\al_{im}+\sum_ma_{itm}\b
x_{it}\tr\b\be_m+\sum_h\sum_ma_{i,t-1,h}a_{itm}\ga_{hm})}{\sum_{\bl
b_t}\exp(\sum_mb_{tm}\al_{im}+\sum_mb_{tm}\b
x_{it}\tr\b\be_m+\sum_h\sum_ma_{i,t-1,h}b_{tm}\ga_{hm})},\quad
t=1,\ldots,T,
\]
where the sums $\sum_h$ and $\sum_m$ are extended to $1,\ldots,M-1$
and $\b b_t$ is an $(M-1)$-dimensional binary vector with elements
$b_{tm}$. This vector has $M$ possible configurations, corresponding
to the possible configurations of any $\b a_{it}$. Then, the
conditional distribution of $\b y_i$, given $\al_i,\b X_i$ and
$y_{i0}$, is equal to
\begin{equation}
p(\b y_i|\al_i,\b X_i,y_{i0}) =
\frac{\exp(\sum_ma_{i+m}\al_{im}+\sum_t\sum_ma_{itm}\b
x_{it}\tr\b\be_m+\sum_h\sum_ma_{i\times
hm}\ga_{hm})}{\prod_t\sum_{\bl
b_t}\exp(\sum_mb_{tm}\al_{im}+\sum_mb_{tm}\b
x_{it}\tr\b\be_m+\sum_h\sum_ma_{i,t-1,h}b_{tm}\ga_{hm})},
\label{p-joint-multi}
\end{equation}
with $a_{i+m}=\sum_t a_{itm}$ and $a_{i\times hm}=\sum_t
a_{i,t-1,h}a_{itm}$.

Proceeding along the same lines as in Section 3.1, we have to
approximate the logarithm of the denominator of
(\ref{p-joint-multi}) through a first-order Taylor expansion
around $\al_i=0$, $\b\be=\b 0$ and $\ga=0$. We have that
\begin{eqnarray*}
&&\sum_t\log[\sum_{\bl
b_t}\exp(\sum_mb_{tm}\al_{im}+\sum_mb_{tm}\b
x_{it}\tr\b\be_m+\sum_h\sum_ma_{i,t-1,h}b_{tm}\ga_{hm})]\approx\\
&&\hspace*{1cm}\approx\sum_t\bigg[\log(M)+\frac{1}{M}\sum_m(\al_{im}+\b
x_{it}\tr\b\be_m)\bigg]+\frac{1}{M}\sum_ma_{i*m}\ga_{m+},
\end{eqnarray*}
with $a_{i*m}=\sum_ta_{i,t-1,m}$ and $\ga_{m+}$ defined in an obvious
way. Thus the approximating model is
\begin{eqnarray*}
&&p^*(\b y_i|\al_i,\b
X_i,y_{i0})=\\
&&\hspace*{0.5cm}=\frac{\exp(\sum_ma_{i+m}\al_{im}+\sum_t\sum_ma_{itm}\b
x_{it}\tr\b\be_m+\sum_h\sum_ma_{i\times
hm}\ga_{hm}-\sum_ma_{i*m}\ga_{m+}/M)}{\sum_{\bl
B}\exp(\sum_mb_{+m}\al_{im}+\sum_t\sum_mb_{tm}\b
x_{it}\tr\b\be_m+\sum_h\sum_mb_{\times
hm}\ga_{hm}-\sum_mb_{*m}\ga_{m+}/M)},
\end{eqnarray*}
where the sum at the denominator is extended to all the possible
configurations of the binary matrix $\b B = \pmatrix{\b b_1 & \cdots
& \b b_T}$ and $b_{+m}$, $b_{\times hm}$ and $b_{*m}$ are defined in
an obvious way.

It may be easily realized that $a_{i+m}$ are sufficient statistics
for the incidental parameters $\al_{im}$ ($i=1,\ldots,n$,
$m=1,\ldots,M-1$) and so, as usual, we can rely on the conditional
distribution
\[
\frac{\exp(\sum_t\sum_ma_{itm}\b
x_{it}\tr\b\be_m+\sum_h\sum_ma_{i\times
hm}\ga_{hm}-\sum_ma_{i*m}\ga_{m+}/M)}{\sum_{\bl
B}^*\exp(\sum_t\sum_mb_{tm}\b
x_{it}\tr\b\be_m+\sum_h\sum_mb_{\times
hm}\ga_{m+}-\sum_mb_{*m}\ga_{m+}/M)},
\]
to estimate the structural parameters, where the sum $\sum^*_{\bl B}$
is extended to al the matrices $\b B$ such that $b_{+m}=a_{i+m}$,
$m=1,\ldots,M-1$.

\section{Conclusions}
We proposed an estimation approach for dynamic logit models for
binary panel data allowing for unobserved heterogeneity and lagged
response variable beyond strictly exogenous covariates. The approach
is based on approximating the assumed logit model with a quadratic
exponential model (Cox, 1972). On the basis of the latter we
construct an approximate conditional likelihood which does not depend
on the heterogeneity parameters, which are considered as incidental
parameters. By maximizing this likelihood, we obtain an approximate
conditional estimator for the other parameters of the logit model,
i.e. the parameters for the covariates and that for the state
dependence, which are referred to as structural parameters. We also
show how this estimator may be improved by using a more precise
approximation of the assumed logit model. The resulting estimator is
the one we suggest to use in practical applications.

The main feature of the estimator above is that it is simpler to use
and performs better than other conditional estimators existing in the
literature. In particular, with respect to the weighted conditional
estimator of Honor\'{e} \& Kyriazidou (2000), that we consider a
benchmark estimator in this literature, our estimator does not
require a kernel function for weighting the response configurations,
may also be used when $T\geq 2$, instead of $T\geq 3$, and in the
presence of time dummies, without requiring particular adjustments. A
more important aspect is that, usually, our estimator also has a
smaller bias and a greater efficiency. This conclusion is based on a
simulation study that we performed along the same lines as Honor\'{e}
\& Kyriazidou (2000). In particular, we noticed that our estimator
has always a limited bias. It also has a root mean square error and a
median absolute error that decrease, as $n$ grows, at a rate close to
$\sqrt{n}$. Moreover, the advantage in terms of bias and efficiency
over the estimator of Honor\'{e} \& Kyriazidou (2000) is more
consistent when there is a strong state dependence effect. An
intuitive explanation of the better performance of our estimator over
their estimator is that the first is based on a conditional
likelihood to which a larger number of response configurations
contribute (actual sample size) with respect to the likelihood on
which the other estimator is based. The larger actual sample size
more than compensate the fact that we are relying on an approximate
conditional likelihood.

In our approach, we also show how it is possible to estimate standard
errors for the proposed estimator. These standard errors are
estimated in the usual way on the basis of an information matrix
which is obtained as a by-product from the estimation algorithm. On
the basis of these standard errors we can construct confidence
intervals for the structural parameters. As our simulation study
shows, these confidence intervals usually have an actual coverage
level very close to the nominal one and so we conclude that the
suggested method for estimating the standard errors is adequate in
practical applications. For this reason, we had not the exigence to
develop more sophisticated methods, based for instance on a bootstrap
procedure, for estimating the standard errors.

In the present paper, we also outlined the extension of the approach
to more complex structures for the state dependence, based on more
than one lagged response variables among the regressors, and to that
of dynamic multinomial logit models for categorical response
variables having more than two categories. We reserve the development
of both of them and the assessment of the quality of the inference
produced in these cases to future research.
\section*{Appendix}
\noindent{\sc Proof of Theorem 1.}
First of all consider that,
under the quadratic exponential model (\ref{quad}), we can express
the conditional distribution of any $\b y_i$, given $\al_i$, $\b
X_i$ and $y_{i0}$, as
\[
p^*(\b y_i|\al_i,\b X_i,y_{i0}) = \frac{\exp(-0.5y_{i0}\ga)}{\mu_{it}}\prod_t
\eta_{it}(y_{i,t-1},y_{it}),
\]
with
$\eta_{it}(y_{i,t-1},y_{it})=\de_{it}(y_{it})\exp(y_{i,t-1}y_{it}\ga)$
and $\de_{it}(y_{it})=\exp(y_{it}\al_i+y_{it}\b
x_{it}\tr\b\be-0.5y_{it}\ga)$ if $t<T$ and
$\de_{it}(y_{it})=\exp(y_{it}\al_i+y_{it}\b x_{it}\tr\b\be)$ if
$t=T$. Therefore, by marginalizing with respect to any response
variable in backward order (from $t=T$), we obtain
\[
p^*(y_i,\ldots,y_{it}|\al_i,\b X_i,y_{i0}) =
\frac{\exp(-0.5y_{i0}\ga)}{\mu_{it}}\left[\prod_{s\leq t}
\eta_{is}(y_{i,s-1},y_{is})\right]g_{i,t+1}(y_{it}),\quad t=1,\ldots,T-1,
\]
where, since $\eta_{it}(y_{i,t-1},0)$ is always equal to $1$, the function
$g_{it}(y_{i,t-1})$ is defined recursively as
\[
g_{it}(y_{i,t-1}) =\left\{\begin{array}{ll}
1+\eta_{iT}(y_{i,T-1},1) & \mbox{if } t = T\\
g_{i,t+1}(0)+\eta_{it}(y_{i,t-1},1)g_{i,t+1}(1) & \mbox{if } t<T
\end{array} \right..
\]
We therefore have that
\[
\frac{p^*(y_{i1},\ldots,y_{it}|\al_i,\b
X_i,y_{i0})}{p^*(y_{i1},\ldots,y_{i,t-1}|\al_i,\b X_i,y_{i0})}
=\frac{\left[\prod_{s\leq t}
\eta_{is}(y_{i,s-1},y_{is})\right]g_{i,t+1}(y_{it})}{\left[\prod_{s\leq t-1}
\eta_{is}(y_{i,s-1},y_{is})\right]g_{it}(y_{i,t-1})}=\eta_{it}(y_{i,t-1},y_{it})
\frac{g_{i,t+1}(y_{it})}{g_{it}(y_{i,t-1})},
\]
which does not depend on $y_{i0},\ldots,y_{i,t-2}$ and so $y_{it}$ is
conditional independent on these variables given $\al_i$, $\b X_i$, $y_{i0}$
and $y_{i,t-1}$. From this conditional probability, expression (\ref{logit})
directly follows.

Finally, on the basis of a Taylor series expansion around $\al_i=0$, $\b\be=\b
0$ and $\ga=0$, we obtain
\[
\log[g_{iT}(y_{i,T-1})]\approx\log(2)+0.5(\al_i+\b x_{iT}\tr\b\be+y_{i,T-1}\ga)
\]
and then
\[
g_{iT}(y_{i,T-1})\approx 2\exp[0.5(\al_i+\b
x_{iT}\tr\b\be)]\exp(0.5y_{i,T-1}\ga)=\exp(c_{iT})\exp(0.5y_{i,T-1}\ga),
\]
with $c_{iT}$ denoting a constant term with respect to $y_{i,T-1}$.
By substituting the latter in $g_{i,T-1}(y_{i,T-2})$ and following
the same recursion above with the Taylor approximation used at any
iteration, we obtain
\[
g_{it}(y_{i,t-1})\approx \exp(c_{it})\exp(0.5y_{i,t-1}\ga),\quad t=1,\ldots,T.
\]
The approximation $\log[g_{it}(1)/g_{it}(0)]\approx 0.5\ga$ then
follows.\vspace*{0.5cm}

\noindent{\sc Proof of Theorem 2}: Let $\widehat{Q}_n(\b\th) =
\ell^*(\b\th)/n$ and $Q_0(\b\th)=E_0\{\log[p^*(\b y|\al,\b
X,y_0,y_+)]\}$. We first prove existence and consistency of
$\hat{\b\th}$ and then asymptotic normality.
\begin{itemize}
\item ({\em Existence and consistency}) Under our assumptions,
conditions (i), (ii) and (iii) of Theorem 2.7 of Newey \& McFadden
(1994) are satisfied and then, since $\hat{\b\th}_n = {\rm
argmax}_{\bl\th} \widehat{Q}_n(\b\th)$, we have that $\hat{\b\th}_n$
exists with probability 1 as $n\rightarrow\infty$ and
$\hat{\b\th}_n\stackrel{p}{\rightarrow}\b\theta_0$. In particular:
\begin{enumerate}
\item[(i)] {\em $Q_0(\b\theta)$ is uniquely maximized at
$\b\theta_0$.} Using a notation derived from Section 4.1, let $\b u(\b D,y_0,\b
y)=(\sum_{t>1}y_t\b d_t\tr, -0.5y_*+y_\times)\tr$. The first derivative of
$Q_0(\b\th)$ at $\b\th_0$ may be then expressed as
\begin{equation}
\b\nabla_{\bl\th} Q_0(\b\th_0)=E_0\{\b u(\b D,y_0,\b y)-E_0[\b u(\b D,y_0,\b
y)|\b D,y_0,y_+]\}=\b 0.\label{der0}
\end{equation}
Moreover, the second derivative may be expressed as
\begin{equation}
\b\nabla_{\bl\th\bl\th} Q_0(\b\th_0)=-E_0[\b S(\b D,y_0,y_+)],\label{der2}
\end{equation}
where $\b S(\b D,y_0,y_+)=V_0[\b u(\b D,y_0,\b y)|\b D,y_0,y_+]$,
with $V_0(\cdot)$ denoting the variance-covariance operator under the
true model. Note however that $\b u(\b D,y_0,\b y)$ may also be
expressed as $\b A(\b D)\b w(y_0,\b y)$, with
\[
\b A(\b D) = \pmatrix{\b D & \b 0\cr \b 0\tr  & 1}\quad\mbox{and}\quad \b
w(y_0,\b y) = \pmatrix{\b y_{-1}\cr -0.5y_*+y_{\times}}
\]
and $\b y_{-1}$ denoting the reduced vector $\b y$ without the first
element. Therefore, (\ref{der2}) may also be expressed as $-E_0\{\b
A(\b D)V_0[\b w(y_0,\b y)|\b D,y_0,y_+]\b A(\b D)\tr\}$, which exists
and is negative definite provided that $E_0(\b D\b D\tr)$ exists and
is of full rank. This is because $V_0[\b w(y_0,\b y)|\b D,y_0,y_+]$
is positive definite for any $y_0$ and $\b D$ and any $y_+$ between 0
and $T$, but the probability that $0<y_+<T$ is always positive.
%FB: bisogna ricorcardi nel frattempo di chiedere a un matematico
\item[(ii)] {\em $\b\th_0$ is an element of the interior of a convex
set $\b\Th$ and $\widehat{Q}_n(\b\th)$ is concave}. That $\b\th_0$ is
an interior point of $\b\Th$ obvious since $\b\Th=\mathbb{R}^{k+1}$.
The concavity of $\widehat{Q}_n(\b\th)$ directly derives from the
concavity of $\ell^*(\b\th)$ discussed at the end of Section 4.1.
\item[(iii)] {\em $\widehat{Q}_n(\b\th)\stackrel{p}{\rightarrow}
Q_0(\b\th)$ for any $\b\th\in\b\Th$}. Since $\widehat{Q}_n(\b\th)$ is the
sample mean of random variables, each with the same expected value equal to
$Q_0(\b\th)$, this easily follows from the law of large number. Note, in
particular, that this law may be applied since $Q_0(\b\th)$ exists for any
$\b\th$ which, in turns, directly derives from the existence of $E_0[\b u(\b
D,y_0,y_+)]$ ensured by that of $E_0(\b D\b D\tr)$.
\end{enumerate}
\item ({\em Normality}) It follows form Theorem 3.1 of Newey \&
McFadden (1994). In particular, the following conditions of this
Theorem hold:
\begin{enumerate}
\item[(i)] {\em $\hat{\b\th}_n\stackrel{p}{\rightarrow}\b\th_0$
and $\b\th_0$ belongs to the interior of $\b\Th$} (see the proof above).
\item[(ii)] {\em $\widehat{Q}_n(\b\th)$ is twice continuously
differentiable in a neighborhood $\cg N$ of $\b\th_0$}. This derivative is
equal to  minus the information matrix (\ref{inf}) divided by $n$ which is
clearly continuous in any $\cg N$.
\item[(iii)] {\em
$\sqrt{n}\,\b\nabla_{\bl\th}\widehat{Q}_n(\b\th_0)\stackrel{d}{\rightarrow}N(\b
0,\b\Sigma)$.} First of all we have that, because of (\ref{der0}),
$E_0[\b\nabla_{\bl\th}\widehat{Q}_n(\b\th_0)]=\b 0$. This implies that
$V_0[\b\nabla_{\bl\th}\widehat{Q}_n(\b\th_0)]=
E_0\{\b\nabla_{\bl\th}\widehat{Q}_n(\b\th_0)\b\nabla_{\bl\th}\widehat{Q}_n(\b\th_0)\tr\}$.
The latter may however be expressed as
\[
E_0\{E_0[\b\nabla_{\bl\th}\widehat{Q}_n(\b\th_0)\b\nabla_{\bl\th}\widehat{Q}_n(\b\th_0)\tr|\b
D,y_0,y_+]\} = E_0\{V_0[\b u(\b D,y_0,\b y)|\b D,y_0,y_+]\},
\]
which, in turn, is equal the $\b\Sigma=-\b\nabla_{\bl\th\bl\th}
Q_0(\b\th_0)$ which exists and is positive definite. The
convergence to the Normal distribution therefore follows from the
Central Limit Theorem.
\item[(iv)] {\em $\sup_{\bl\th\in\cgl
N}\|\b\nabla_{\bl\th\bl\th}\widehat{Q}_n(\b\th)+\b\Sigma\|\stackrel{p}{\rightarrow}0$.}
This directly follows from Lemma 2.4 of Newey \& McFadden (1994) and
the fact that $E_0[\b\nabla_{\bl\th\bl\th}\widehat{Q}_n(\b\th_0)] =
-\b\Sigma$ and that
$E_0[\|\b\nabla_{\bl\th\bl\th}\widehat{Q}_n(\b\th)\|]$ is finite for
any $\b\th\in\cg N$.
\item[(v)] {\em $\b\Sigma$ is nonsingular}. See item (iii) above.
\end{enumerate}
\end{itemize}
\section*{References}
{\small\begin{description}
\item  Agresti, A. (2002): {\em Categorical data analysis}. New
York: John Wiley \& Sons.\vspace*{-0.2cm}
\item Arellano, M. and Honor\'e, B. (2001): ``Panel Data Models:
Some Recent Developments'', in {\em Handbook of Econometrics},
Vol. V,
 Ed. Heckman J. J. and  Leamer E.. Amsterdam: North-Holland.\vspace*{-0.2cm}
\item Andersen, E. B. (1970): ``Asymptotic properties of
conditional maximum-likelihood estimators'', {\em Journal of Royal
Statistical Society, B}, 32, 283--301.\vspace*{-0.2cm}
\item Andersen, E. B.  (1972):  ``The numerical solution of a set
of conditional estimation equations'', {\em Journal of the Royal
Statistical Society,  B}, 34, 42--54.\vspace*{-0.2cm}
\item Chamberlain, G. (1985): ``Heterogeneity, omitted variable bias,
and duration dependence'', in  {\em Longitudinal analysis of labor
market data}, Ed. Heckman J. J. and  Singer B.. Cambridge: Cambridge
University Press.\vspace*{-0.2cm}
%
%\item Chamberlain, G., (1993): ``Feedback in panel data models'',
%unpublished manuscript, Department of Economics, Harvard
%University.\vspace*{-0.2cm}
%
\item Cox, D. R. (1972): ``The analysis of multivariate binary
data'', {\em Applied Statistics}, 21, 113--120.\vspace*{-0.2cm}
\item Cox, D. R. and  Wermuth, N. (1994): ``A note on the
quadratic exponential binary distribution'', {\em Biometrika}, 81,
403--408.\vspace*{-0.2cm}
\item Heckman, J. J. (1981a): ``Statistical models for discrete
panel data'', in {\em Structural Analysis of Discrete Data}, Ed.
McFadden D. L. and Manski C. A..  Cambridge, MA: MIT
Press.\vspace*{-0.2cm}
\item Heckman, J. J. (1981b): ``Heterogeneity and state
dependence'', in {\em Structural Analysis of Discrete Data}, Ed.
McFadden D. L. and Manski C. A.. Cambridge, MA: MIT
Press.\vspace*{-0.2cm}
%//
\item Honor\'e, B. E. and Kyriazidou, E. (2000): ``Panel data
discrete choice models with lagged dependent variables'', {\em
Econometrica}, 68, 839--874.\vspace*{-0.2cm}
\item Hyslop, D. R. (1999): ``State dependence, serial correlation
and heterogeneity in intertemporal labor force participation of
married women'', {\em Econometrica},  67,
1255--1294.\vspace*{-0.2cm}
\item Hsiao, C. (1986): {\em Analysis of Panel Data}. New York:
Cambridge University Press.\vspace*{-0.2cm}
\item Manski, C. (1987):  ``Semiparametric analysis of random
effects linear models from binary panel data'', {\em
Econometrica}, 55, 357--362.\vspace*{-0.2cm}
\item Newey W. K. and McFadden D. L. (1994): ``Large Sample
Estimation and Hypothesis Testing'', in  {\em Handbook of
Econometrics},  Vol. 4,  Ed. Engle R. F. and McFadden D. L..
Amsterdam: North-Holland.\vspace*{-0.2cm}
\item Neyman, J. and Scott, E. L. (1948): ``Consistent estimates
based on partially consistent observations'', {\em Econometrica},
16, 1--32.\vspace*{-0.2cm}
\item Rasch, G. (1961): ``On general laws and the meaning of
measurement in psychology'', {\em Proceedings of the IV Berkeley
Symposium on Mathematical Statistics and Probability}, 4,
321--333.\vspace*{-0.2cm}
\end{description}}

\end{document}